\documentclass[article]{siamart220329}

\usepackage{amssymb}
\usepackage{graphicx}

\usepackage{hyperref}    
\usepackage{subcaption}
\usepackage{amssymb,latexsym,amsmath}
\usepackage{xcolor}
\usepackage[%
    style=numeric-comp,maxbibnames=99,sorting=nyt,
    sortcites=true,doi=false,url=true,   giveninits=true,hyperref]{biblatex}

\usepackage{tikz,tikz-cd}
\usetikzlibrary{hobby}
\usetikzlibrary{arrows}
\usetikzlibrary{positioning}
\usetikzlibrary{shapes,snakes}
\usetikzlibrary{cd}

\def\R{\mathbb{R}}
\def\bone{{\bf 1}}
\def\xvec{{\bf x}}

\def\N{\mathbb{N}}
\DeclareMathOperator{\diag}{diag}
\DeclareMathOperator{\dd}{dd}

\newtheorem{remark}[theorem]{{\sc Remark}}

\numberwithin{equation}{section}

\newcommand\newnew[1]{{\color{black}#1}}

\begin{document}

\title{Weighted enumeration of nonbacktracking walks on weighted graphs \thanks{Submitted to the editors DATE.}}

\author{
Francesca Arrigo\thanks{Department of Mathematics and Statistics, University of Strathclyde, Glasgow, UK, G1 1XH (\email{francesca.arrigo@strath.ac.uk}). The work of F.A. was supported by fellowship ECF-2018-453 from the Leverhulme Trust.}
\and Desmond J. Higham \thanks{School of Mathematics, University of Edinburgh, James Clerk Maxwell Building, Edinburgh, UK, EH9 3FD (\email{d.j.higham@ed.ac.uk}). The work of D.J.H. was supported 
the Engineering and Physical Sciences Research Council 
under grants EP/P020720/1 and EP/V015605/1.}
\and Vanni Noferini\thanks{Aalto University, Department of Mathematics and Systems Analysis, P.O. Box 11100, FI-00076, Aalto, Finland (\email{vanni.noferini@aalto.fi}). Supported by an Academy of Finland grant (Suomen Akatemian p\"{a}\"{a}t\"{o}s 331230).}
\and Ryan Wood\thanks{Corresponding author. Aalto University, Department of Mathematics and Systems Analysis, P.O. Box 11100, FI-00076, Aalto, Finland (\email{ryan.wood@aalto.fi}). Supported by an Academy of Finland grant (Suomen Akatemian p\"{a}\"{a}t\"{o}s 331240)}}

\date{}

\maketitle
\begin{abstract} 
    We extend the notion of nonbacktracking walks from unweighted graphs to graphs whose edges have a nonnegative weight. 
    Here the weight associated with a walk is taken to be the product over the weights along the individual edges. 
    We give two ways to compute the associated generating function, and corresponding node centrality measures. One 
    method works directly on the original graph and one uses a line graph construction
    followed by a projection. The first method is more efficient, but the second has the advantage of extending naturally to 
    time-evolving graphs. \newnew{Based on these generating functions, we define and study corresponding centrality measures.} Illustrative computational results are also provided.
\end{abstract}

\begin{keywords}
Complex network, matrix function, generating function, line graph, combinatorics, evolving graph, temporal network, centrality measure, Katz centrality.
\end{keywords}

\begin{MSCcodes}
 05C50, 05C82, 68R10
\end{MSCcodes}

\section{Introduction}\label{sec:intro}
Complex network analysis is an expanding scientific discipline that has recently been producing many research challenges, with applications across
several fields of science and engineering \cite{Estradabook,Newmanbook}. One important question is that of ranking the nodes of a graph by importance, or, in more mathematical terms, defining and studying an appropriate \emph{centrality measure}. Those centrality measures that can be formulated and computed via combinatorial properties of walks on the underlying graph have received special attention  \cite{BK13,EV05,fenu2017block,Katz53,MRMPO10} because they have convenient formulations in terms of linear algebra
that lead to efficient computational methods. In recent years, this paradigm has been refined by studying centrality measures that are based on counting not all walks but only some of them, namely, walks that do not backtrack \cite{AGHN17a,AGHN17b,Beyond,GHN18,TCDM21,TCTE21,TSE19a} or more generally do not cycle \cite{Beyond}. Nonbacktracking walks are known to be linked to zeta functions of graphs \cite{HST,KKS21,MS01,ST96}. Their associated centrality measures have been shown to possess attractive computational properties 
\cite{AGHN17a,AGHN17b,AHN19b,GHN18}
and have been studied both for undirected and directed graphs, and more recently for time evolving graphs \cite{AHNW22}. However, in the context of \newnew{the combinatorics of} nonbacktracking walks, so far only unweighted graphs have been studied. \newnew{We mention that nonbacktracking random walks were previously considered in \cite{Kempton16}, but the problems studied there are different to the ones analyzed in the present paper. Moreover, \cite{Kempton16} focuses only on the special case where the nodes are given themselves a positive weight $\varphi(i)$, and the weight of the edge $(i,j)$ is defined as $\omega((i,j))=\varphi(i)\varphi(j)$. Instead, we do not impose any restriction on the edges' weights.} In the theory of graph zeta functions, weighted graphs have been considered  by defining the weight of a walk to be the sum (and not the product, as in this paper) of the weights of its edges \cite{HST}. We discuss this issue further in section~\ref{sec:background}.

The main purpose of the present paper is to extend the combinatorial theory of nonbacktracking walks, 
\newnew{and corresponding centrality measures}, to graphs whose edges carry a positive \emph{weight}. These graphs are associated with generic nonnegative adjacency matrices, in contrast to unweighted graphs that correspond to binary adjacency matrices. 
While for unweighted graphs one may be interested in the enumeration of walks of a given length, for weighted graphs the combinatorial problem is more sophisticated due to the presence of weights. The edge weights naturally give rise to 
an overall weight for each walk, a concept that can be used 
alongside the length (i.e., the number of edges traversed). 

The structure of the paper is as follows. 
In Section~\ref{sec:background} we introduce some relevant notation and core concepts.
Section~\ref{sec:genfun} sets up and studies the 
issue of characterizing the classical generating function 
associated with nonbacktracking weighted walks
\newnew{and using it to compute a centrality measure}.
In section~\ref{sec:linegraph}
we introduce an alternative formulation 
that applies to a wider class of 
generating functions \newnew{and centrality measures}.
Section~\ref{sec:evolve} shows how these ideas can be extended to the case of evolving graphs. 
\newnew{Numerical experiments are conducted in Section~\ref{sec:num}.}
We finish in Section~\ref{sec:disc} with a brief discussion. 

\section{Background and Notation}\label{sec:background}
In this paper, we consider finite graphs. A finite graph is a triple $G=(V,E,\Omega)$ where $V=[n]$ is the set of the nodes (or vertices), $E \subset V \times V$ is the set of (directed) edges, and $\Omega : E \rightarrow (0,\infty)$ is a weight function that associates to each edge a positive weight. If $\Omega(e)=1$ for all $e \in E$, then the graph is said to be {\it unweighted}; if for any pair $i\neq j$, with $i,j \in V$ we have that $(i,j) \in E \Leftrightarrow (j,i) \in E$ and that $\Omega((i,j))=\Omega((j,i))$ then the graph is said to be {\it undirected}; and if, for every $i \in V$, we have that $(i,i) \not\in E$ then the graph is said to be without loops. Graphs that are not unweighted are usually called {\it weighted} and graphs that are not undirected are usually called {\it directed}. It is, however, convenient (and we will do so within this work) to relax the terminology so that the set of directed (resp., weighted) graphs contains as special cases also undirected (resp., unweighted) graphs, further we will assume that all graphs are without loops. 

A \textit{walk of length $\ell$} on the graph $G$ is a sequence of nodes $i_1,i_2, \dots, i_{\ell+1}$ such that $(i_j,i_{j+1})\in E$ for all $1 \leq j \leq \ell$. Equivalently, it can be seen as a sequence of edges $e_1,\dots,e_\ell$ such that  $e_j \in E$ for all $j=1,\dots,\ell-1$ and the end node of $e_j$ coincides with the starting node of $e_{j+1}$. 

\begin{definition}\label{def:weightofwalk}
Let $G=(V,E,\Omega)$ be a weighted graph. 
The {\rm weight} of the walk $e_1,\dots,e_\ell$ is 
\[ \prod_{k=1}^\ell \Omega(e_k)\]
where $\Omega(e_k)$ is the weight of the edge $e_k\in E$.
\end{definition}
\begin{remark}
When $\Omega: E\to \{1\}$ is the weight function associated with an unweighted graph, then the weight of all walks in the network is one, regardless of their length.
\end{remark}
In the context of mainstream graph theory, the \textit{weight} (or \textit{length} or \textit{cost}) of a walk is sometimes defined as the \emph{sum}, rather than the product, of the weights of its edges. In that scenario, zeta functions of graphs (which are closely related to the enumeration of nonbacktracking walks) have been studied \cite{HST}. However, we argue that within complex network analysis 
Definition~\ref{def:weightofwalk} \newnew{has several useful applications}. For example, consider a road network where nodes represent towns and a nonnegative integer 
edge weight 
$A_{ij}$ records the number of distinct roads
connecting town $i$ and town $j$.
Then, the number of distinct routes from $i$ to $j$ that pass through one intermediate town is equal to
\[
\sum_{k=1}^n A_{ik}A_{kj},
\]
that is, the weighted sum of walks of length two,  where the weight is the \emph{product} of the weights of its edges. Similarly, in a model where edges represent independent probabilistic events and their weights are their probabilities, as discussed in the original work of Katz \cite{Katz53}, it is natural to postulate that the weight of a walk is the product of the weight of its edges, in agreement with the fact that the joint probability of a sequence of independent events is the product of the individual probabilities.

Given a node ordering, the corresponding 
adjacency matrix of a graph is the matrix 
$A \in \R^{n \times n}$ entrywise defined as: 
\[ A_{ij} = \begin{cases}
0 \ &\mathrm{if} \ (i,j) \not \in E;\\
\Omega((i,j)) \ &\mathrm{if} \ (i,j) \in E.
\end{cases}\]
Note that a graph is undirected if and only if its adjacency matrix is symmetric; it is without loops if and only if its adjacency matrix has zero diagonal; and it is unweighted if and only if its adjacency matrix has entries all lying in $\{0,1\}$.

The problems of \newnew{enumerating walks in} unweighted graphs and \newnew {enumerating} \newnew{weighted} walks in weighted graphs may both be solved by considering powers of the adjacency matrix: indeed, the $(i,j)$ entry of $A^k$ is equal to, respectively, the number of walks of length $k$ from node $i$ to node $j$ (when the graph is unweighted)  or the weighted sum of walks of length $k$ from node $i$ to node $j$ (when the graph is weighted). As a consequence, the generating function for the (possibly weighted) enumeration of walks is given by
\[ I + t A + t^2 A^2 + \dots = \sum_{k=0}^\infty t^k A^k = (I-tA)^{-1} ,\]
where we adopt the standard convention that the (weighted) sum of walks of length zero from $i$ to $j$ is 1 if $i=j$ and 0 otherwise.
Here, $t$ is a real parameter small enough to ensure convergence of the series which scales by $t^k$ the count for walks of length $k$.  

A walk can also be seen as a sequence of nodes. If the sequence does not contain a subsequence of the form $iji$ for some nodes $i$ and $j$, then the walk is said to be \emph{nonbacktracking} (NBT). 
\newnew{We} define $p_k(A)$ \newnew{to be} the matrix whose $(i,j)$ entry contains the sum of the weights of all nonbacktracking walks of length $k$ from node $i$ to node $j$.
By convention, $p_0(A) = I$.
Note that, by definition, $p_k(A) \leq A^k$ elementwise. 
Combinatorially, the problem of computing the (weighted) enumeration of nonbacktracking walks is equivalent to finding an explicit expression for the generating function
\begin{equation}\label{eq:gggoal}
\Phi(t) = 
     \sum_{k=0}^\infty t^k p_k(A)
\end{equation}
for suitable values of the parameter $t>0$.

This problem was addressed in \cite{GHN18} for unweighted undirected graphs, and later in \cite{AGHN17a} for unweighted directed graphs. In \cite{AGHN17b}, the solution was extended to the more general generating function 
\begin{equation}\label{eq:gggoal2}
     \kappa(t) = \sum_{k=0}^\infty c_k t^k p_k(A),
\end{equation}
where $(c_k)_k \subset [0,\infty)$ is an arbitrary sequence. In \cite{AHNW22}, the theory was further extended to consider time evolving graphs. However, so far, the 
quantities \eqref{eq:gggoal} and \eqref{eq:gggoal2} have not yet been studied for weighted graphs. Their characterization is the main contribution of this paper.

A corollary of obtaining such computable expressions 
is a numerical recipe for associated nonbacktracking centralities. Indeed, beyond its algebraic interpretation as a generating function, \eqref{eq:gggoal2} can be interpreted analytically as a function that will converge for sufficiently small values of the variable $t$. Choosing one such value for $t$ allows us to define a centrality measure based on the weighted sum of edges. For example, if $\bone$ is the vector of all ones, then the $i$-th component of the vector
\[ \left(  \sum_{k=0}^\infty c_k t^k p_k(A) \right) \bone \]
computes a nonbacktracking version of Katz centrality \cite{Katz53}. \newnew{The latter is} defined as the doubly weighted sum of all the walks departing from node $i$, where the weight of each walk within the sum is 
the product of the weight of the walk itself and $t^k$, where $k$ is the walk length. Similarly, for the {\it subgraph centrality} version of nonbacktracking Katz, the doubly weighted sum of all the walks that start and end on node $i$ is given by
\[ \left(  \sum_{k=0}^\infty c_k t^k p_k(A) \right)_{ii}. \] 
As a consequence, two additional questions that we address in this paper are to describe the radius of convergence of \eqref{eq:gggoal2} 
and to derive computable expressions 
for the associated centrality measures.
\newnew{We refer to \cite{AGHN17a,AGHN17b,GHN18}, and the references therein, for details of the benefits of
nonbacktracking 
in the centrality context.}

We consider two approaches to bridge the gap between weighted graphs and current results on the combinatorics of nonbacktracking walks. The first is specialized to the case $c_k \equiv 1$, i.e., to compute \eqref{eq:gggoal}; it leads directly 
to an expression that has computational advantages \newnew{as it does not require to go through the edge-level and, thus, it requires the construction of a potentially much smaller matrix than the second approach}. The second is based on a technique, described in \cite{Beyond, AHNW22}, of forming the line graph, obtaining a generating function there, and finally projecting back to compute \eqref{eq:gggoal2}. While, potentially, the second approach may be computationally less efficient, it has the advantages that (i) it is able to solve the more general problem \eqref{eq:gggoal2}, (ii) it can be generalized to the 
setting of time evolving graphs, and (iii) it allows us to easily estimate the convergence radius of \eqref{eq:gggoal2} (including the special case of \eqref{eq:gggoal}).

\section{The generating function of nonbacktracking walks on a weighted graph}\label{sec:genfun}

In this section, we assume that  $G$ is a finite directed weighted graph with $n$ nodes, without loops, and having adjacency matrix $A$. The directed edge from node $i$ to node $j$ has weight $A_{ij} > 0$. 
Following Definition~\ref{def:weightofwalk}, to the walk 
$   i_1 i_2 i_3 \dots i_{\ell+1}  $
of length $\ell$ we assign 
the weight
$ A_{i_1i_2} A_{i_2 i_3} \cdots A_{i_\ell i_{\ell+1}}$. 
We note the distinction here between the \emph{length}
and the \emph{weight} of a walk.

The goal of this section is to obtain a convenient formula for the generating function
$\Phi(t)$ in (\ref{eq:gggoal}).
We note that this generalizes the 
version previously studied for an 
unweighted graph \cite{AGHN17a}, and 
the expression $\Phi(t) \mathbf{1}$
is then a natural candidate for a node centrality measure.

\subsection{Describing the matrices $p_k(A)$ via a recurrence relation}

Let us first set up some further notation: given two square matrices $X,Y\in\mathbb{R}^{n\times n}$ we distinguish between matrix multiplication, $XY$, and elementwise multiplication, $X\circ Y$, where $(X\circ Y)_{ij} = X_{ij}Y_{ij}$. 
Similarly, we differentiate between the $k$-th linear algebraic power $X^k$ and the $k$-th 
elementwise power, $X^{\circ k}$, so $(X^{\circ k})_{ij} = (X_{ij})^k$.
Moreover, following Matlab notation, $\dd(X):=\diag(\diag(X))$ will denote the diagonal matrix whose diagonal entries are equal to the diagonal entries of $X$. We
first prove the following $k$-term recurrence, which 
generalizes previous results that have been derived 
independently for the 
unweighted
\cite{BL70,TP09} and 
undirected \cite{ST96} cases.

\begin{theorem}\label{thm:recurrence}
For all $k \geq 1$,
\[ p_k(A) = \sum_{\substack{\ell=2h+1 \ \mathrm{odd} \\ 1 \leq \ell \leq k}} (A^{\circ (h+1)} \circ (A^T)^{\circ h} ) p_{k-\ell}(A) - \sum_{\substack{\ell=2h \ \mathrm{even} \\ 2 \leq \ell \leq k}} \dd((A^{\circ h})^2) p_{k-\ell}(A).  \]
\end{theorem}
\begin{proof}
For the base case of $k=1$, the statement reduces to $p_1(A)=(A \circ \bone \bone^T)p_0(A)$. Since $A\circ \bone \bone^T=A$ and $p_0(A)=I$, in turn this yields $p_1(A)=A$, which is manifestly true since any walk of length one is nonbacktracking. Let us now give a proof by induction.

We start by considering $A p_{k-1}(A)$, whose $(i,j)$ entry is equal to the sum of the weights of all walks of length $k$ from $i$ to $j$ that are nonbacktracking if the first step is removed. This value is equal to $p_k(A)_{ij}$ plus the sum of the weights of all backtracking walks of length $k$ from $i$ to $j$ that are nonbacktracking if the first step is removed. Such walks must be of the form $iai\dots j$: the weight of one such walk is $A_{ia} A_{ai}$ times the weight of a certain NBT walk of length $k-2$ from $i$ to $j$. Summing over all $a$ adjacent to $i$ yields $\dd(A^2) p_{k-2}(A)$. However, we have subtracted too much, because any such walk of the form $iaia \dots j$, being backtracking after removing the first step, was not present in $(A p_{k-1}(A) )_{ij}$. The weight of one such walk is $A_{ia}A_{ai}A_{ia} = A_{ia}^2 A_{ai}$ times the weight of a certain NBT walk of length $k-3$ from $a$ to $j$. We can sum again over all $a$ adjacent to $i$, to obtain $((A^{\circ 2} \circ A^T) p_{k-3}(A))_{ij}.$ We should sum this value back, but again we are adding a bit too much, because walks satisfying the previous requirements and being of the form $iaiai \dots j$ should not be there.

It is clear that this sequence of corrections goes on until we exhaust the length of the walk and the statement of the theorem is a consequence of the two following facts, both true for all $h \geq 0$.
\begin{enumerate}
\item The total weight of walks of length $k$ from $i$ to $j$ of the form $i(ai)^{h}a  \dots j$, such that the final subwalk (of length $k-(2h+1)$) from $a$ to $j$ is not backtracking, is equal to
\[ \sum_{a : (i,a) \in E} (A_{ia})^{h+1} (A_{ai})^h p_{k-2h-1}(A)_{aj} = \left((A^{\circ (h+1)} \circ (A^T)^{\circ h}) p_{k-2h-1}(A)\right)_{ij}. \]
\item The total weight of walks of length $k$ from $i$ to $j$ of the form $(ia)^{2h} i  \dots j$, such that the final subwalk (of length $k-2h$) from $i$ to $j$ is not backtracking, is equal to
\[ \sum_{a : (i,a) \in E} (A_{ia})^{h}(A_{ai})^h p_{k-2h}(A)_{ij} = \left(\dd((A^{\circ h})^2) p_{k-2h-1}(A)\right)_{ij}. \]
\end{enumerate}         
\end{proof}

\subsection{Solving the recurrence relation}

Let us continue by giving a combinatorial result \newnew{in Proposition \ref{prop:growrec}. Its statement} expresses the generating function of a sequence satisfying a growing recurrence relation \newnew{in terms of two individual generating functions}.
\begin{proposition}\label{prop:growrec}
Let $(\mathcal{P}_k)_k$ and $(\mathcal{C}_\ell)_\ell$ be two sequences in some (possibly noncommutative) ring, and suppose that $(\mathcal{P}_k)_k$ satisfies the growing recurrence
\[ \sum_{\ell=0}^k \mathcal{C}_\ell \mathcal{P}_{k-\ell} = 0 \]
for all $k \geq 1$. Then, the (formal) generating functions $\Phi(t)=\sum_{k=0}^\infty \mathcal{P}_k t^k$ and $\Psi(t)=\sum_{\ell=0}^\infty \mathcal{C}_\ell t^\ell$ are related by the formula $\Psi(t)\Phi(t)=\mathcal{C}_0\mathcal{P}_0.$ 
\end{proposition}
\begin{proof}
Observe that, using the recurrence,
\[ \Psi(t) \Phi(t) = \sum_{k=0}^\infty t^k \sum_{\ell=0}^\infty \mathcal{C}_\ell \mathcal{P}_{k-\ell} = \mathcal{C}_0 \mathcal{P}_0. \] 
\end{proof}

We can now apply the general technique of Proposition~\ref{prop:growrec} to the special case of the generating function \eqref{eq:gggoal}, whose coefficients satisfy the recurrence described in Theorem~\ref{thm:recurrence}. 
In other words, we specialize Proposition~\ref{prop:growrec} to sequences in the ring $\mathbb{R}^{n \times n}$ where $\mathcal{P}_k = p_k(A)$ and
\[ \mathcal{C}_\ell = \begin{cases} I \ &\mathrm{if} \ \ell=0;\\
-[A^{\circ(h+1)} \circ (A^T)^{\circ h}] \ &\mathrm{if} \ \ell=2h+1;\\
\dd((A^{\circ h})^2) \ &\mathrm{if} \ \ell=2h>0.
\end{cases}\]
In particular, $\mathcal{C}_0=\mathcal{P}_0=I$, and hence by Proposition \ref{prop:growrec} $\Phi(t)=\Psi(t)^{-1}$.
In turn, we can write $\Psi(t)=\Psi_e(t)-\Psi_o(t)$ by splitting even and odd terms and by extracting the minus sign appearing in the odd terms of $(\mathcal{C}_\ell)_\ell$. It is easy to see that $\Psi_e(t)$ is diagonal while $\Psi_o(t)$ is the off-diagonal part of $\Psi(t)$, since we assume $G$ to be without loops. 
Moreover, 
\[ \left(\Psi_o(t)\right)_{ij} = \sum_{h=0}^\infty t^{2h+1} A_{ij}^{h+1} A_{ji}^h = \frac{t A_{ij}}{1-t^2 A_{ij} A_{ji}}.  \]
Similarly,
\[ \left(\Psi_e(t)\right)_{ii} = 1 + \sum_{h=1}^\infty t^{2h} \sum_{j=1}^n A_{ij}^h A_{ji}^h =     1 + \sum_{j=1}^n \frac{t^2 A_{ij} A_{ji}}{1-t^2 A_{ij} A_{ji}}.\]
Let $S=A\circ A^T$,  \newnew{let $Q=S^{\circ 1/2}$}, let
\[ 
f_1(x)=\frac{x}{1-x}, \quad f_2(x)=\frac{x}{1+x},
\]
and let $f_i(\circ tX)$ denote the \emph{elementwise} application of $f_i$ to the matrix $t X$, for $i=1,2$. 
Then, if we denote by $\circ/$ the elementwise application of $/$, we can write
\[ \Psi_e(t) = I + \dd(f_1(\circ t Q) f_2 (\circ t Q)), \qquad \Psi_o(t) = t A \circ/ (\bone \bone^T-t^2 S) \]
and hence $\Psi(t) = I + \dd(f_1(\circ t Q) f_2 (\circ t Q)) - t A \circ/ (\bone \bone^T-t^2 S)$. 

We can state this more formally as a theorem.
\begin{theorem}\label{thm:genkatz}
\label{thm:Phi}
In the notation above, 
for all values of $t$ such that \eqref{eq:gggoal} converges, we have 
\begin{equation}\label{eq:goalachieved}
\Phi(t) = (I + \dd[f_1(\circ t Q) f_2 (\circ t Q)] - t A \circ/ (\bone \bone^T-t^2 S) )^{-1}.
\end{equation} 
\end{theorem}

As a sanity check, let us see what happens in three distinct interesting limiting cases that have been addressed previously in the literature.
\begin{itemize}
\item  First, let us verify that in the limit of an unweighted graph we recover \cite[equation~(3.3)]{AGHN17a}. In this case, $(A_{ij})^h=A_{ij} \in \{0,1\}$ for all $h \geq 1$. As a result, if $D=\dd(A^2)$,
\[ \left(\Psi_e(t)\right)_{ii} = 1 +\sum_{h=1}^\infty t^{2h} \sum_{j=1}^n A_{ij} A_{ji} = 1 + D_{ii} \frac{t^2}{1-t^2} \Rightarrow \Psi_e(t) = \frac{I-t^2 I + t^2 D}{1-t^2} \]
and
 \[ \left(\Psi_o(t)\right)_{ij} =  t A_{ij} + \sum_{h=1}^\infty t^{2h+1} A_{ij} A_{ji} = t A_{ij} + \frac{t^3 S_{ij}}{1-t^2} \Rightarrow \Psi_o(t) = \frac{t A - t^3 (A-S)}{1-t^2}  \]
 which imply the known $\Phi(t) = (1-t^2) (I - t A + t^2(D-I) + t^3 (A-S))^{-1}$ from \cite[Equation~(3.3)]{AGHN17a}.
 \item Next, let us observe that if no edge is reciprocated, that is, if there is no $(i,j)\in E$ such that $(j,i)\in E$, then $S=Q=0$. Hence, we recover the generating function associated with classical Katz centrality, i.e., $\Phi(t)=(I-tA)^{-1}$, which is consistent
 with the fact that every walk is nonbacktracking under this assumption.
 \item  Finally, if the graph is undirected then $S=A^{\circ 2}$ and $Q=A$. Hence, the formulae simplify to
\[ \Psi_e(t) = I + \dd[f_1(\circ t A) f_2 (\circ t A)], \qquad \Psi_o(t) = t A \circ/ (\bone \bone^T-t^2 A^{\circ 2})\]
yielding in particular
\[ \Phi(t)=(I + \dd[f_1(\circ t A) f_2 (\circ t A)] - t A \circ/ (\bone \bone^T-t^2 A^{\circ 2}))^{-1}. \]
If we additionally assume that the graph is unweighted, we further reduce to $\Phi(t)=(1-t^2)(I-At+t^2(D-I))^{-1}$ in agreement with \cite[Equation~(5.3)]{GHN18}.
\end{itemize}

We now briefly comment on the convergence of {$\Phi(t) = \sum_k p_k(A) t^k$} to the right-hand-side of \eqref{eq:goalachieved}. 
Since the series converges to a rational function, its radius of convergence is equal to the smallest of its poles. One way to compute the radius is therefore via the eigenvalues of the rational function ${\Psi(t)=}\Phi(t)^{-1}$. A more straightforward method (albeit possibly less efficient) to {estimate} the radius of convergence is available when computing $\Phi(t)$ with a different method. \newnew{This is described in more detail in Section~\ref{sec:genfun} and, in particular, within Corollary~\ref{cor:radius}}. In spite of the somewhat awkward notation, \eqref{eq:goalachieved} is in fact quite straightforward to compute given $A$, by composing elementwise functions and matrix addition and multiplications. 


We conclude this section by recalling that we can define a nonbacktracking version of Katz centrality on weighted graphs by summing the value of the generating function over all possible ending nodes, which can be expressed as the linear algebraic matrix-vector multiplication $\Phi(t) \mathbf{1}$.

The following corollary is then an immediate consequence of  
Theorem~\ref{thm:recurrence}.

\begin{corollary}
\label{cor:cent}
For all values of $t$ such that \eqref{eq:gggoal} converges,
consider the centrality measure where node $i$ is assigned the value 
$x_i$ according to $\xvec{} = \Phi(1) \mathbf{1}$.
Then $\xvec$ may be found by solving the linear system
\begin{equation}
 (I + \dd[f_1(\circ t Q) f_2 (\circ t Q)] - t A \circ/ (\bone \bone^T-t^2 S) )
 \xvec{} 
 = \mathbf{1}.
 \label{eq:xcent}
\end{equation}
\end{corollary}

Corollary~\ref{cor:cent} shows in particular that the centrality measure
can be found without explicitly computing the inverse in \eqref{eq:goalachieved}. We can instead compute the vector of nonbacktracking centralities $\bf x$ by solving the linear system 
(\ref{eq:xcent}).
We note that the coefficient matrix in 
(\ref{eq:xcent}) is 
no less sparse than $I-tA$; hence the computational complexity of solving such a linear system is the same as for classical Katz centrality, and 
the task is feasible \newnew{with standard tools for sparse linear systems} for very large, sparse networks.

\section{Generating function by constructing the line graph and projecting back}\label{sec:linegraph}

In this section, we derive an alternative computable expression 
for the generating function $\Phi(t)$. Although generally this second method is less computationally efficient, it offers three
main advantages: (i) it can be extended to nonbacktracking centrality measures 
other than Katz 
(for example, based the exponential rather than the resolvent); 
(ii) it allows for a simple characterization of the radius of convergence of the generating function; and (iii) it can be extended to time evolving graphs.

As before, we consider a finite weighted graph with $n$ nodes. We also assume (directed) edges have been labelled from 
$1$ to $m$ in an arbitrary, but fixed, manner. 
We may then define the \textit{source matrix} $L \in \R^{m \times n}$ and \textit{target} (or \textit{terminal}) \textit{matrix} $R \in \R^{m \times n}$ as follows~\cite{VDF09}:
$$ L_{ej}=\begin{cases}
1 \ \ \mathrm{if} \ \mathrm{edge} \ e \ \mathrm{starts} \ \mathrm{from} \ \mathrm{node} \ j\\
0 \ \ \mathrm{otherwise}
\end{cases} \qquad R_{ej}=\begin{cases}
1 \ \ \mathrm{if} \ \mathrm{edge} \ e \ \mathrm{ends} \ \mathrm{on} \ \mathrm{node} \ j\\
0 \ \ \mathrm{otherwise}.
\end{cases} $$
Moreover, we let \newnew{$Z$} be an $m\times m$ diagonal matrix such that \newnew{$Z_{ee} = A_{ij}$}, where (in the chosen labelling of the edges) the $e$-th edge is precisely $(i,j)$.\footnote{For clarity, we will sometimes use the notation $i\to j$ to denote the edge $(i,j)\in E$.} Then, we have the following relationship.  
\begin{proposition}\label{prop1}
We have \newnew{$A = L^T Z R$}.
\end{proposition}
\begin{proof}
Since $Z$ is diagonal, $(L^T Z R)_{ij} = \sum_{e=1}^m L_{e i} Z_{ee} R_{e j}$. But there is at most one value of $e$ such that $L_{ei}R_{ej} \neq 0$, and that is precisely the value identifying the edge $i \rightarrow j$, if this is an edge of the graph. If such an edge does not exist then the summation yields $0$, as desired. If such an edge exists, then, for that $e$, $Z_{ee}=A_{ij}$ which concludes the proof.
\end{proof}

Now let $W$ be the weighted matrix of the dual graph (or line graph), i.e., the graph whose nodes correspond to the original (directed) edges, and whose edges are pairs of edges from the original graph that can form a walk. The pair $(i \rightarrow j, j \rightarrow k)$ represents a walk that has weight equal to the product of the original edge weights; that is, $A_{ij} A_{jk}$. These values are recorded in the entries of $W$, with $W_{ef} = A_{ij}A_{jk}$ if $e$ is the label of edge $i\to j$ and $f$ is the label of edge $j\to k$.
\begin{theorem}\label{thm0}
We have \newnew{$W = Z R L^T Z.$}
\end{theorem}
\begin{proof}
We proceed entrywise. 
Suppose for concreteness that edge $e$ is $i \rightarrow j$ and edge $f$ is $k \rightarrow \ell$, where $i \neq j$, $k \neq \ell$ are (possibly, but not necessarily, all distinct) nodes. 
Note for a start that $W_{ef} = A_{ij}A_{j\ell}$ if $j=k$ and $W_{ef}=0$ if $j \neq k$. 
Now, since $Z$ is diagonal,
$$(Z R L^T Z)_{ef} = Z_{ee} Z_{ff} \sum_{h=1}^n R_{e h} L_{f h} = A_{ij} A_{k\ell} \sum_{h=1}^n R_{eh} L_{fh}.$$
Suppose $j \neq k$; then there is no $h$ such that $R_{eh}L_{fh} \neq 0$, so the summation above is $0=W_{ef}$. On the other hand, if $j=k$ then the summation over $h$ yields $1$ so that $(ZRL^TZ)_{ef} =A_{ij}A_{j\ell} = W_{ef}$.
\end{proof}

In the unweighted case, we have a projection relation $L^T W^k R=A^{k+1}$ \cite[Proposition~2.4]{Beyond}. However, for weighted graphs, entries of  $W^k$ count walks of length $k+1$, but with incorrect weights. For example, the walk $1 \rightarrow 2 \rightarrow 3 \rightarrow 4$ would be weighted $A_{12} A_{23}^2 A_{34}$ rather than $A_{12} A_{23} A_{34}$. We now exhibit a trick that corrects this problem. Coherently with the notation of the previous section, below $M^{\circ 1/2}$ denotes the elementwise nonnegative square root of a nonnegative matrix $M$; note that generally this does not correspond to the classical matrix square root $\sqrt{M}$ (i.e., the matrix $X$ such that $X^2=M$), a notable exception being the case of a diagonal square matrix with nonnegative diagonal. We note that $Z$ falls in this latter  category, hence the notation in the following result. 

\begin{theorem}\label{thm1}
Let $0<k  \in \N$. The $(e,f)$ element of \newnew{$\sqrt{Z} (W^{\circ 1/2})^k \sqrt{Z}$} counts, with weights, all walks of length $k+1$ from edge $e$ to edge $f$.
\end{theorem}
\begin{proof}
 The crucial observation  is that $W^{\circ 1/2} = \sqrt{Z} R L^T \sqrt{Z}$, which is clear by a minor modification of the proof of Theorem \ref{thm0}. We now proceed by induction on $k$. For the base case $k=1$, it suffices to observe that $\sqrt{Z} W^{\circ 1/2} \sqrt{Z} = Z R L^T Z = W.$
Suppose now that the statement holds for $k-1$. Then,
$$ \sqrt{Z} (W^{\circ 1/2})^k \sqrt{Z} = \sqrt{Z} (W^{\circ 1/2})^{k-1} \sqrt{Z} (Z^{-1/2}) W^{\circ 1/2} \sqrt{Z}.$$
Define for notational simplicity $U: =   \sqrt{Z} (W^{\circ 1/2})^{k-1} \sqrt{Z}$, $X:= Z^{-1/2}$, $Y:=W^{\circ 1/2}$, $\Sigma:=\sqrt{Z}$. Then, since $X$ and $\Sigma$ are diagonal,
$$(UXY\Sigma)_{ef} = \sum_{g \in E} U_{e g} X_{gg} Y_{gf} \Sigma_{ff}.$$
Suppose now that edge $e$ is $i \rightarrow j$ and edge $f$ is $h \rightarrow \ell$; then edges $g$ must be of the form $x \rightarrow h$ for some node $x$. Indeed, $Y_{gf} = 0$ unless the end node of edge $g$ coincides with the 
start node of edge $f$, i.e., unless $gf$ is a walk of length two. Hence, in this notation,
$$(UXY\Sigma)_{ef} = \sum_{x:A_{xh}>0}U_{e, x\to h}\sqrt{A_{h\ell}A_{xh}} \sqrt{\frac{A_{h\ell}}{A_{xh}}} = A_{h\ell}\sum_{x: A_{xh}>0} U_{e, x\to h},$$
where 
$\sum_{x: A_{xh} > 0} U_{e, x\to h}$ is, by the inductive assumption, the count (with weights) of all walks of length $k-1$ from edge $e$ to all edges of the form $x \rightarrow h$, i.e., the weighted enumeration of all walks of length $k-1$ from edge $e$ to node $h$. However, the count with weights of all walks of length $k$ from edge $e$ to edge $f$ is precisely the count with weights of all walks of length $k$ from edge $e$ to node $\ell$ with node $h$ as the penultimate node, i.e., the right hand side in the latter displayed equation.
\end{proof}

We have the following consequence of Theorem \ref{thm1}. 
\begin{corollary}\label{cor:Akp1}
For all $k \in \N$, \newnew{$L^T\sqrt{Z} (W^{\circ 1/2})^k \sqrt{Z} R = A^{k+1}$.}
\end{corollary}
\begin{proof}
The result follows from Proposition~\ref{prop1}, if $k=0$, and from Theorem~\ref{thm1}, if $k>0$.
\end{proof}

Now let $B\in\mathbb{R}^{m\times m}$ be the nonbacktracking version of $W$, i.e., $B_{ef} = 0$ if $W_{ef} W_{fe} \neq 0$ and $B_{ef} = W_{ef}$ otherwise. This matrix is often referred to as the \textit{Hashimoto matrix}~\cite{Hash90}. 
Recall, moreover, that $p_{k}(A)\in\mathbb{R}^{n\times n}$ is the matrix counting all NBT walks of length $k$ (from $i$ to $j$ in its $(i,j)$ element). Now we can observe that all the proofs above hold for $B$ as well, modulo substituting walks with nonbacktracking walks. Hence, the projection relation still holds.
\begin{theorem}\label{thm2}
For all $k \in \N$, 
we have 
\newnew{$L^T \sqrt{Z} (B^{\circ 1/2})^k \sqrt{Z} R = p_{k+1}(A).$}
\end{theorem}
\begin{proof}
We have 
$p_1(A)=A=L^T Z R$, and when $k>0$ the result follows from a 
minor modification of the arguments used to prove 
Corollary~\ref{cor:Akp1}.
\end{proof}

Suppose now that $(c_k)_k \subset [0,\infty)$ is a sequence and $t$ is such that 
$$\kappa(t,A) = \sum_{k=0}^\infty c_k t^kp_{k}(A)$$
as in \eqref{eq:gggoal2} converges; we are interested in the centrality measure
\begin{equation}
\boldsymbol{v}(t,A) = \kappa(t,A)\mathbf{1}.
\label{eq:nbtcen}
\end{equation}
\newnew{We now derive formulae for $\kappa(t,A)$ and $\boldsymbol{v}(t,A)$. 
To this end, we introduce the following notation. Given a real-analytic scalar function
$$f(x) = \sum_{k=0}^\infty c_k x^k$$
consider the 
operator
$$\partial f(x) = \sum_{k=0}^\infty c_{k+1} x^k = \frac{f(x) - c_0}{x}.$$
We then have the following.
\begin{theorem} \label{thm:line}
It holds that 
\[ \sum_{k=0}^\infty c_k t^k p_k(A) = c_0 I + t L^T \sqrt{Z} \partial f(t V) \sqrt{Z} R, \]
for $V=B^{\circ 1/2}$ and $|t|<r/\rho(V)$, where $\rho(V)$ is the spectral radius of $V$ and $r$ is the radius of convergence of the scalar function $f(x)=\sum_{k=0}^\infty c_k t^k$.

Hence, for the centrality associated with $f(x)$ and $t$ small enough to give convergence in the matrix power series, in \eqref{eq:nbtcen} we have 
$$\boldsymbol{v}(t,A) = c_0\mathbf{1} +t  L^T \sqrt{Z} \partial f(t V)   \sqrt{Z} \mathbf{1}.$$
\end{theorem}
\begin{proof}
    By Theorem~\ref{thm2} we easily see that
$$ \kappa(t,A)   = c_0 I +t  L^T \sqrt{Z}\left( \sum_{k=0}^\infty c_{k+1} t^k (B^{\circ 1/2})^k\right) \sqrt{Z} R.$$

    As a consequence,
$$\boldsymbol{v}(t,A) = c_0\mathbf{1} +t  L^T \sqrt{Z}\left( \sum_{k=0}^\infty c_{k+1} t^k (B^{\circ 1/2})^k\right) \sqrt{Z} \mathbf{1}.$$
\end{proof}
}

Observing that the resolvent is an eigenfunction (with eigenvalue $1$) of $\partial$, we note in particular that for Katz centrality, i.e., $c_k = 1$ for all $k$, $\partial f(x) = f(x) = (1-x)^{-1}.$
Hence, we have the following special case.
\begin{corollary}\label{cor:Phigen}
\newnew{In the notation of Theorem~\ref{thm:line}, we have that} the generating function $\Phi(t)$ defined in \eqref{eq:gggoal} can be expressed as
\begin{equation}\label{eq:goalachieved2}
\newnew{\Phi(t) = I + t L^T \sqrt{Z} (I-tV)^{-1} \sqrt{Z} R.}
\end{equation}
\end{corollary}

This analysis in particular yields a lower bound for the radius of convergence for \eqref{eq:gggoal}.

\begin{corollary}\label{cor:radius}
If $|t|<\rho(V)^{-1}$, where $V=B^{\circ 1/2}$, then the sequence $\Phi(t)=\sum_{k=0}^\infty p_k(A) t^k$ converges.
\end{corollary}

\newnew{
\begin{remark}
    Letting $r$ denote the radius of convergence of \eqref{eq:gggoal}, Corollary~\ref{cor:radius} shows that $r \geq \rho(V)^{-1}$. It is possible to strengthen this result and prove that $r = \rho(V)^{-1}$. A proof of this fact, which is beyond the scope of the present article, 
    appears in \cite[Theorem 5.2]{manmevanni}.
\end{remark}}

\section{Nonbacktracking centralities for evolving weighted graphs}
\label{sec:evolve} 

In Sections~\ref{sec:genfun} and \ref{sec:linegraph}, we obtained formulae for the generating function $\Phi(t)$ in \eqref{eq:gggoal} by working, respectively, at node and edge level. For a static network, i.e., one which does not evolve in time,  working at the node level is clearly preferable as, for large $n$, we may have that $n \ll m$. However, a significant advantage of the latter, edge-level, formula is that it easily extends to the case of temporal networks in all backtracking regimes, whereas a direct node-level formula which forbids backtracking in time is generally unavailable \cite{AHNW22}. Let us first 
generalize the definition of graph, walk, and NBT walk to the dynamic case.

\begin{definition}\label{def:timeevolving}
A finite time-evolving graph $\mathcal{G}$ is a finite collection of graphs  $ (G^{[1]}, \dots, G^{[N]})$, associated with the non-decreasing time stamps $(t_1, \dots, t_N) \in \mathbb{R}^{N}$, such that the set of nodes of $G^{[i]}$ does not depend on $i$ and when observed at time $t_i$ the structure of $\mathcal{G}$ is identical to that of $G^{[i]}$.
\end{definition}

We remark that the concept of a graph can be extended to 
the dynamic setting in a number of ways \cite{holme11}.
The discrete-time framework of Definition~\ref{def:timeevolving} covers a 
range of realistic scenarios where interactions take place, or are recorded, at specific points in time. For example, 
in an on-line social media platform, an edge may represent a 
form of communication between users, and 
$G^{[i]}$ may count the number of interactions between each pair of individuals over time $(t_{i-1},t_i]$.

The definition of walk across a network can be extended to the setting of temporal graphs as follows. 
\begin{definition}\label{def:Twalk1}
{\rm A walk of length $\ell$} across a temporal network is defined as an ordered sequence of $\ell$ edges $e_1 e_2 \dots e_{\ell}$ such that for all $k =1,\ldots,\ell-1$ the end node of $e_k$ coincides with the start node of $e_{k+1}$  and, moreover, that $e_k \in E^{[\tau_1]}, e_{k+1} \in E^{[\tau_2]}$ for some $1\leq \tau_1 \leq \tau_2\leq N$, where $E^{[\tau_i]}$ denotes the set of edges of the graph $G^{[\tau_i]}$.
\end{definition}
It is useful to make an equivalent definition.
\begin{definition}\label{def:Twalk2}
{\rm A walk of length $\ell$} across a temporal network is defined as an ordered sequence of $\ell+1$ nodes $i_1 i_2 \dots i_{\ell+1}$ such that for all $k = 2,\ldots,\ell$ it holds that $i_{k-1}\to i_{k} \in E^{[\tau_1]}$ and $i_k\to i_{k+1} \in E^{[\tau_2]}$ for some $1\leq\tau_1 \leq \tau_2\leq N$.
\end{definition}
We want to stress that multiple edges can be crossed at one given time stamp and, moreover, that a walk is allowed to remain inactive   
for some of the time stamps.
We also recall here that there is not just one definition of backtracking for temporal networks; indeed, three arise naturally~\cite{AHNW22}:
\begin{itemize}
\item backtracking happens within a certain time-stamp; we will refer to this as {\it backtracking in space},
\item backtracking happens across time-stamps; we will refer to this as {\it backtracking in time},
\item backtracking happens both within a time-stamp and across time-stamps (not necessarily in this order); we will refer to this as {\it backtracking in time and space}.
\end{itemize}

Given any finite time-evolving graph $\mathcal{G}$, we can associate with it a matrix $M$ called the \textit{global temporal transition matrix} which was defined in~\cite{AHNW22} for unweighted graphs. Definition~\ref{def:M} below generalizes the definition of the global temporal transition matrix to the weighted case.

\begin{definition}\label{def:M}
Let ${\mathcal{G}} = (G^{[1]},G^{[2]},\ldots,G^{[N]})$ be a time-evolving graph with $N$ time stamps. 
The {\rm weighted global temporal transition matrix} associated with ${\mathcal{G}}$ is the $m\times m$ block matrix 
\begin{equation}
M = M^{[1,\dots,N]} = 
\begin{bmatrix}   
{C^{[1]}} & {C^{[1,2]}} & {C^{[1,3]}} & \dots & {C^{[1,N]}} \\ 
0 & {C^{[2]}} & {C^{[2,3]}} & \dots & {C^{[2,N]}} \\
\vdots & & \ddots &  & \vdots \\
\vdots & & & \ddots & \vdots\\
0 & \dots & \dots & 0 & {C^{[N]}}
\end{bmatrix}^{\circ 1/2}, 
\end{equation}
where the definition of the blocks depends on the chosen backtracking regime in the following way:
\begin{itemize}
\item [(i)] $C^{[\tau_1]} = W^{[\tau_1]}$ and $C^{[\tau_1, \tau_2]} = W^{[\tau_1, \tau_2]} := S^{[\tau_1]} R^{[\tau_1]}(L^{[\tau_2]})^T S^{[\tau_2]}$  for all $\tau_1, \tau_2 = 1,2,\dots,N$ ($\tau_1 < \tau_2$) if backtracking in both space and time is permitted;
\item[(ii)] $C^{[\tau_1]} = B^{[\tau_1]}$ and $C^{[\tau_1,\tau_2]} = W^{[\tau_1,\tau_2]}$ for all $\tau_1,\tau_2=1,2,\ldots,N$ ($\tau_1< \tau_2$) if backtracking in space is forbidden but backtracking in time is permitted;
\item[(iii)] $C^{[\tau_1]} = W^{[\tau_1]}$ and $C^{[\tau_1,\tau_2]} = B^{[\tau_1,\tau_2]} := W^{[\tau_1, \tau_2]} -  (W^{[\tau_1, \tau_2]} \circ  {W^{[\tau_2, \tau_1]}}^T)^{\circ 1/2}$ for all $\tau_1,\tau_2=1,2,\ldots,N$ ($\tau_1< \tau_2$) if backtracking in time is forbidden but backtracking in space is permitted; and 
\item[(iv)] $C^{[\tau_1]} = B^{[\tau_1]}$ and $C^{[\tau_1,\tau_2]} = B^{[\tau_1,\tau_2]}  $ for all $\tau_1,\tau_2=1,2,\ldots,N$ ($\tau_1< \tau_2$) if backtracking in both time and space is forbidden.
\end{itemize} 

\end{definition}

It was further shown in \cite{AHNW22} that the global temporal transition matrix provides an accurate way of counting walks in all backtracking regimes across a finite unweighted time-evolving graph and thereby allows for the computation of the (nonbacktracking) Katz centrality $\boldsymbol{v}$ via the formula:
\begin{equation}\label{eq:Mkatz}
    \boldsymbol{v}(t) = (I +  t\mathcal{L}^T(I-tM)^{-1}\mathcal{R}) \boldsymbol{1},
\end{equation}
where $\mathcal{L}, \mathcal{R}$ are the global source and target matrices respectively as defined in \cite[Definition 4.4]{AHNW22}.

To handle the weighted case, we may extend 
Theorem~\ref{thm2} naturally to the global temporal transition matrix $M$ in the following way.
\begin{theorem}\label{prop:M}
For a finite time-evolving graph with $N$-many time frames, \newnew{let the {\rm global weight matrix} $Z$ be defined block-wise as 
\[Z := Z^{[1,2,\dots, N]} = \diag(Z^{[1]}, Z^{[2]}, \dots, Z^{[N]}),  \]
where $Z^{[\tau_i]}$ is the diagonal matrix associated with time stamp $1 \leq \tau_i \leq N$. For each of these matrices, their diagonal entries are given by $Z^{[\tau_i]}_{ee} = {w_e}^{[\tau_i]}$,} with ${w_e}^{[\tau_i]}$ being the weight of edge $e$ at time stamp $\tau_i$. Further let the backtracking regime be fixed such that the weighted global temporal transition matrix $M$ is fixed. Then, for $0 < k \in \mathbb{N}$, the $(e,f)$-th entry of $\sqrt{Z}M^k \sqrt{Z}$ counts, with weights, all permitted walks of length $k+1$ from edge $e$ to edge $f$ across the time evolving graph given the backtracking regime. 
\end{theorem}

\begin{proof}\label{pf:M}
Suppose the backtracking regime is given such that the structure of $M$ is fixed as specified in Definition~\ref{def:M}. We prove the theorem by induction on the length of permitted walks $k \in \mathbb{N}$. Consider the basis case $k = 1$:
\[
\sqrt{Z} M \sqrt{Z} = 
\sqrt{Z}
\begin{bmatrix}   
{C^{[1]}} & {C^{[1,2]}} & {C^{[1,3]}} & \dots & {C^{[1,N]}} \\ 
0 & {C^{[2]}} & {C^{[2,3]}} & \dots & {C^{[2,N]}} \\
\vdots & & \ddots &  & \vdots \\
\vdots & & & \ddots & \vdots\\
0 & \dots & \dots & 0 & {C^{[N]}}
\end{bmatrix}^{\circ 1/2}
\sqrt{Z}. 
\]

The $(e,f)$-th element of this matrix correctly counts the unique walk of length two with weights from edge $e$ to edge $f$. In the following we will omit the temporal superscript, since the indices $e$ and $f$ also uniquely determine the time frame. With this convention:
\begin{align*}
    \left( \sqrt{Z} M \sqrt{Z} \right)_{ef} &= \sum \limits_{r,s} \sqrt{Z}_{er} M_{rs} \sqrt{Z}_{sf}\\
    &= (\sqrt{Z})_{ee} (\sqrt{Z})_{ff} M_{ef} \\
    &= \begin{cases}
    \sqrt{w_e}\sqrt{w_f} \sqrt{w_ew_f}= w_e w_f & ef \ \text{is a permitted walk of length two}\\
    0 &\text{otherwise}.
    \end{cases}
\end{align*}
where we have used the fact that, if $ef$ is an admissible walk of length two, then $(M)_{ef} = \sqrt{w_ew_f}$.

Suppose now that the result holds for $k-1$ and for brevity denote by $P$ the matrix $\sqrt{Z}M^{k-1}\sqrt{Z} $, which, by the inductive assumption, correctly counts in its entries the number of weighted temporal walks of length $k$. 
Then, using the fact that $Z$ is diagonal, 
\begin{align*}
    (\sqrt{Z}M^k \sqrt{Z})_{ef}  
    &= \sum_{r,s,t} \left(\sqrt{Z}M^{k-1}\sqrt{Z}\right)_{er} (Z^{-1/2})_{rs}M_{st} (\sqrt{Z})_{tf}\\
    &= \sum_{r}  P_{er} (Z^{-1/2})_{rr}M_{rf} (\sqrt{Z})_{ff}\\
    &= \begin{cases}w_f \sum_{r} P_{er}   & e \dots rf \ \text{is a permitted walk of length }k+1\\
    0 & \mathrm{otherwise}. \end{cases}
\end{align*}

By the inductive assumption $P_{er}$ counts the permitted weighted walks of length $k$ from edge $e$ to edge $r$. Therefore the above formula does indeed count the weighted permitted walks beginning with edge $e$ and ending on edge $f$ with length $k+1$ correctly.
\end{proof}

\begin{theorem}
Given global source, target and weight matrices, $\mathcal{L}$, $\mathcal{R}$ and \newnew{$Z$} respectively, we can compute $M$ when backtracking is entirely forbidden in two steps:
\begin{enumerate}
    \item $\widehat{M} = \left(\sqrt{Z} (\mathcal{R} \mathcal{L}^T  -  \mathcal{R}\mathcal{L}^T \circ \mathcal{L}\mathcal{R}^T )\sqrt{Z}\right)$;
    \item Obtain $M$ from $\widehat{M}$ by setting all entries below the block diagonal to $0$.
\end{enumerate}
\end{theorem}
\begin{proof}
    The above theorem is easy to prove by observing that the $(i,j)$-th block of the matrix $\mathcal{R}\mathcal{L}^T$ is equal to $R^{[i]} {L^{[j]}}^T $, whereas the $(i,j)$-th block of the matrix $\mathcal{L}\mathcal{R}^T$ is equal to $L^{[i]}{R^{[j]}}^T = (R^{[j]} {L^{[i]}}^T)^T $; whereupon the $(i,j)$-th block of the matrix $\widehat{M}$ becomes
    \begin{align*}
        \widehat{M}_{ij} = \sqrt{Z^{[i]}} \left(R^{[i]} {L^{[j]}}^T  -R^{[i]} {L^{[j]}}^T  \circ (R^{[j]} {L^{[i]}}^T)^T  \right) \sqrt{Z^{[j]}}. 
    \end{align*}
    The central term here in brackets can be seen as the binarized $B^{[i,j]}$, i.e., $B^{[i,j]}$ where all non-zero weights are uniformly equal to $1$, thus the presence of a non-zero entry $(B^{[i,j]})_{ef}$ simply reflects whether or not the concatenation of edges $e$ and $f$ forms a non-backtracking walk of length two.  Matrix multiplication from the left by  $\sqrt{Z^{[i]}}$ and on the right by $\sqrt{Z^{[j]}}$ then provides the appropriate weighting for the $(e,f)$-th entry, namely $\sqrt{w_e} \sqrt{w_f}$, as required. Finally, the second step of setting all blocks below the block diagonal, i.e., $\widehat{M}_{ij}$ with $i > j$, to zero reflects the requirement that walks may not move back in time.
\end{proof}

We can also compute the $f$-total communicability of the time-evolving graph with weights by using the global temporal transition matrix $M$. 
\begin{theorem}\label{thm:fcenttime}
Given a function $f$ with series expansion $f(t) = \sum_{k = 0}^\infty c_k t^k$ having radius of convergence $r$, and some fixed backtracking regime, the $f$-total communicability $\boldsymbol{v}_f(t)$ of the time-evolving graph $\mathcal{G} = (G^{[1]}, \dots, G^{[N]}) $ with $N$ time stamps is given by the formula:
\begin{equation}\label{eq:fcenttime}
    \newnew{\boldsymbol{v}_f(t) = (c_0 I + t \partial \mathcal{L}^T \sqrt{Z}f(tM)\sqrt{Z} \mathcal{R}) \boldsymbol{1}}
\end{equation}
for $0< |t| < r/\max_{i =1, \dots, N}\{\rho(C^{[\tau_i]})\} $.
\end{theorem}
\begin{proof}\label{pf:fcenttime}
By Proposition \ref{prop:M}, we have that $\sqrt{Z} M^k \sqrt{Z}$ counts with weights all walks of length $k+1$, thus
\[
    \boldsymbol{v}_f(t) = (c_0 I + t \partial \mathcal{L}^T \sqrt{Z}f(tM)\sqrt{Z} \mathcal{R})\boldsymbol{1} 
    = c_0 I + t \sum_{k = 0}^\infty c_{k+1} t^k \mathcal{L}^T\sqrt{Z} M^k \sqrt{Z} \mathcal{R} \boldsymbol{1}.\\
\] 
In the above formula we see the number of walks of length $k+1$ correctly counted with weights that are further weighted by the coefficient $c_{k+1}t^{k+1}$, which is provided by the series expansion of $f(x)$.
\end{proof}

\section{Numerical Experiments}\label{sec:num}

In this section we show how the formulae for nonbacktracking Katz centrality from sections \ref{sec:genfun}, \ref{sec:linegraph} and \ref{sec:evolve} may produce significantly different node-rankings for real-world social networks when compared with Katz centrality which permits backtracking walks. We further examine the effect of weighting edges on the rankings produced by both centrality measures. To this end, we consider the Katz centrality formula~\eqref{eq:Mkatz} as applied to one static network and one temporal network, both derived from the same data set (Fauci's email release \cite{benson2021fauciemail}) \footnote{\newnew{The code used in the following analysis can be found at \url{https://github.com/rwood12347/Weighted-enumeration-of-nonbacktracking-walks-on-weighted-graphs}}}. 
\newnew{The original dataset is a collection of over 3000 pages of emails involving Anthony Fauci and his staff during the COVID-19 pandemic. Data includes sender and receivers (including CC'd) of emails, as well as time stamps of when the emails were sent. Both networks used in the following were presented in \cite{benson2021fauciemail}}.

\subsection{Analysis on Static Networks}

In this section, we analyze a static network produced by \cite{benson2021fauciemail} which is both undirected and weighted. We have an edge $(i,j)$ if there exists an
email which involves both nodes $i$ and $j$ as any combination of sender and recipient (including CC'd recipients). The weight assigned to such an edge, $\Omega((i,j))$, is a positive integer equal to the number of such emails that were sent.

In our analysis, we apply Corollary~\ref{cor:cent} to obtain the NBT Katz centrality vector for our network, which is then contrasted with the classical Katz centrality vector for attenuation factor values $t = 0.5/\rho$ and $t = 0.95/\rho$, where $1/\rho$ is the radius of convergence for the respective centrality measure. In particular $1/\rho$ is equal to $1/\rho(A)$, where $A$ is the adjacency matrix of the graph in the case of classical Katz centrality; whereas $1/\rho  \newnew{=} 1/\rho(B^{\circ 1/2 })$ in the case of weighted nonbacktracking Katz centrality \cite[Theorem 5.2]{manmevanni}, where $B$ is the Hashimoto matrix associated to the graph. These values are given in Table~\ref{tab:staticConv}. We also analyze the binarized graph which is produced from the static graph by setting all edge weights to $1$. In the context of the binarized network  $1/\rho $ equals $1/\rho(A)$ in the case of classical Katz centrality, and $1/\rho(B)$ for nonbacktracking Katz centrality, where $A$ and $B$ are the adjacency and Hashimoto matrices associated with the binarized network, respectively.
\begin{table}[t]
    \centering
        \caption{Static network convergence information.}
    \begin{tabular}{c|c}
    \hline
        Non-binarized graph $\rho(B^{\circ 1/2})$ &   $926.9$\\ 

         Binarized graph $\rho(B)$ &  $48.61$ \\ 
      \text{Non-binarized nonbacktracking permitted range of $t$} & $t \in [0,1.079\cdot 10^{-3})$   \\ 
         \text{Binarized nonbacktracking permitted range of $t$} &  $t \in [0,2.057\cdot 10^{-2})$ \\
         \hline
         \hline
         Non-binarized $\rho(A)$ &     $1038$  \\ 
         Binarized $\rho(A)$ &  $51.26$ \\ 
      \text{Non-binarized backtracking permitted range of $t$} & $t \in [0,9.635\cdot 10^{-4})$   \\ 
         \text{Binarized backtracking permitted range of $t$} &  $t \in [0,1.951\cdot 10^{-2})$ \\
         \hline
    \end{tabular}
    \label{tab:staticConv}
\end{table}


\begin{figure*}[htbp]
        \centering
        \begin{subfigure}[h]{0.475\textwidth}
            \centering
            \includegraphics[width=1\textwidth]{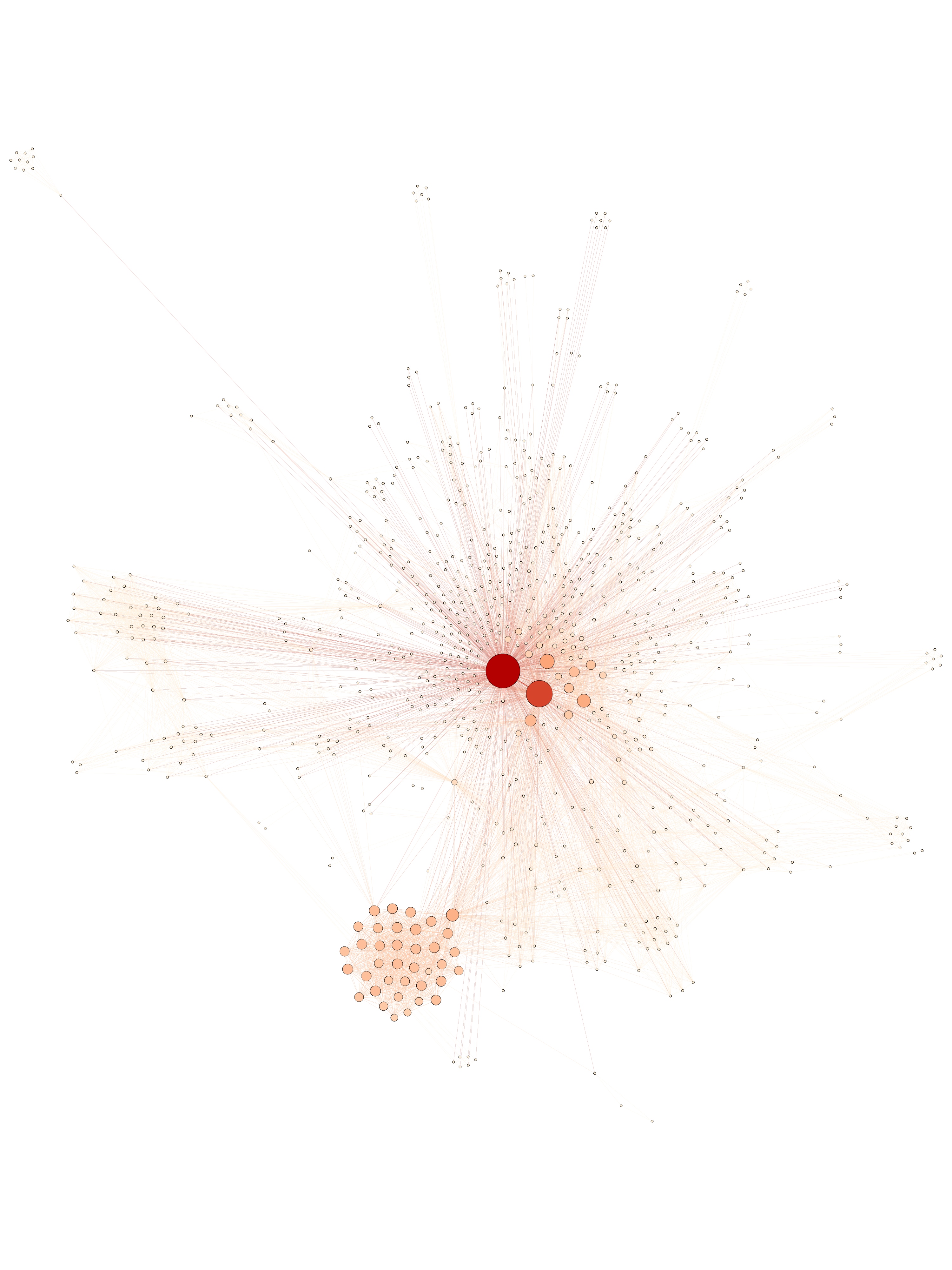}

            {{\small Katz centrality $t=0.5/\rho(A)$}}    
            \label{fig:Katz5}
        \end{subfigure}
        \hfill
        \begin{subfigure}[h]{0.475\textwidth}  
            \centering 
            \includegraphics[width=1\textwidth]{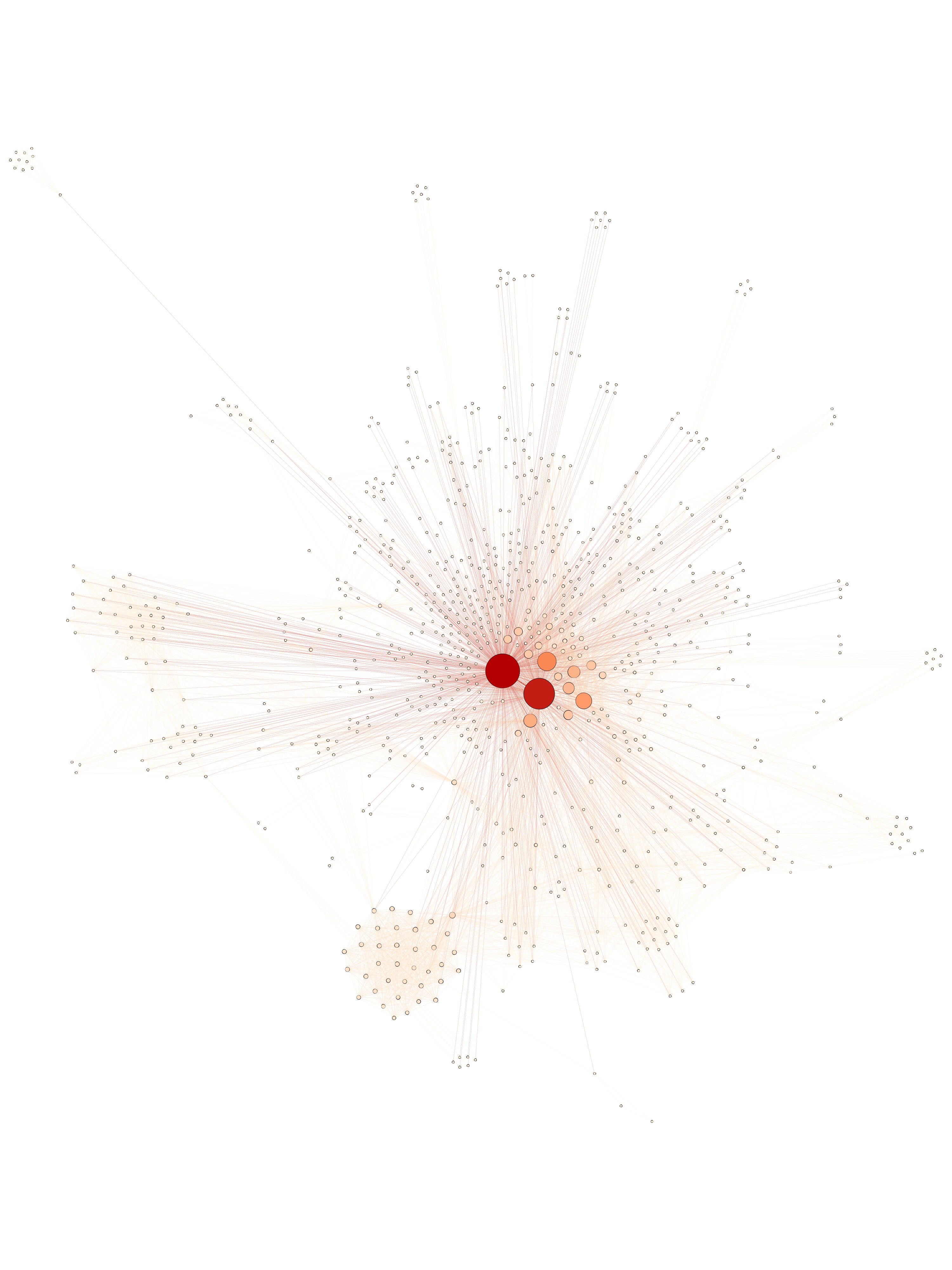}
            
            {{\small Katz centrality $t = 0.95/\rho(A)$}}    
            \label{fig:Katz95}
        \end{subfigure}
        \vskip\baselineskip
        \begin{subfigure}[h]{0.475\textwidth}   
            \centering 
            \includegraphics[width=1\textwidth]{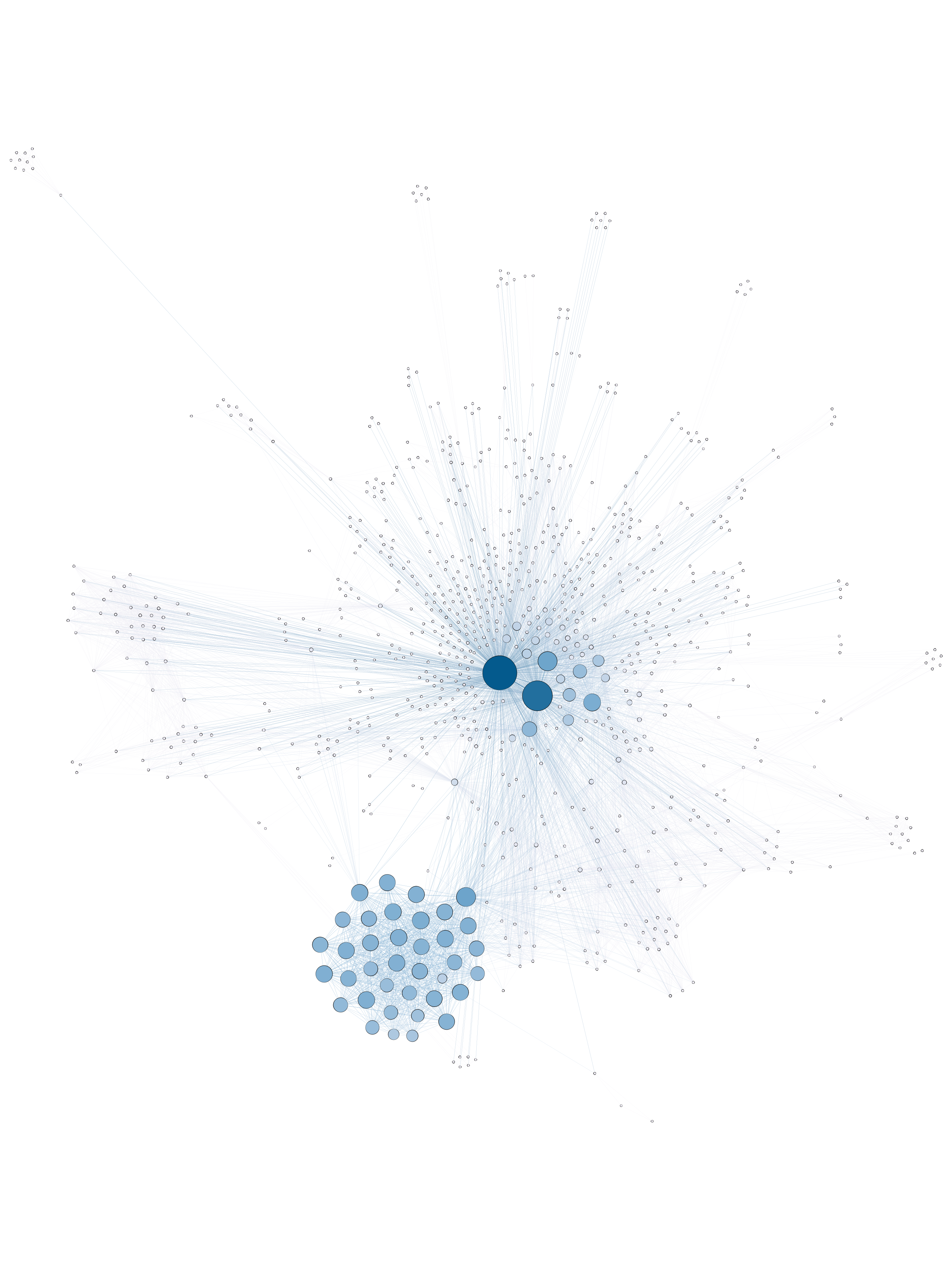}
            
            {{\small NBT Katz centrality $t = 0.5/\rho({B^{\circ 1/2}})$}}    
            \label{fig:NBT5}
        \end{subfigure}
        \hfill
        \begin{subfigure}[h]{0.475\textwidth}   
            \centering 
            \includegraphics[width=1\textwidth]{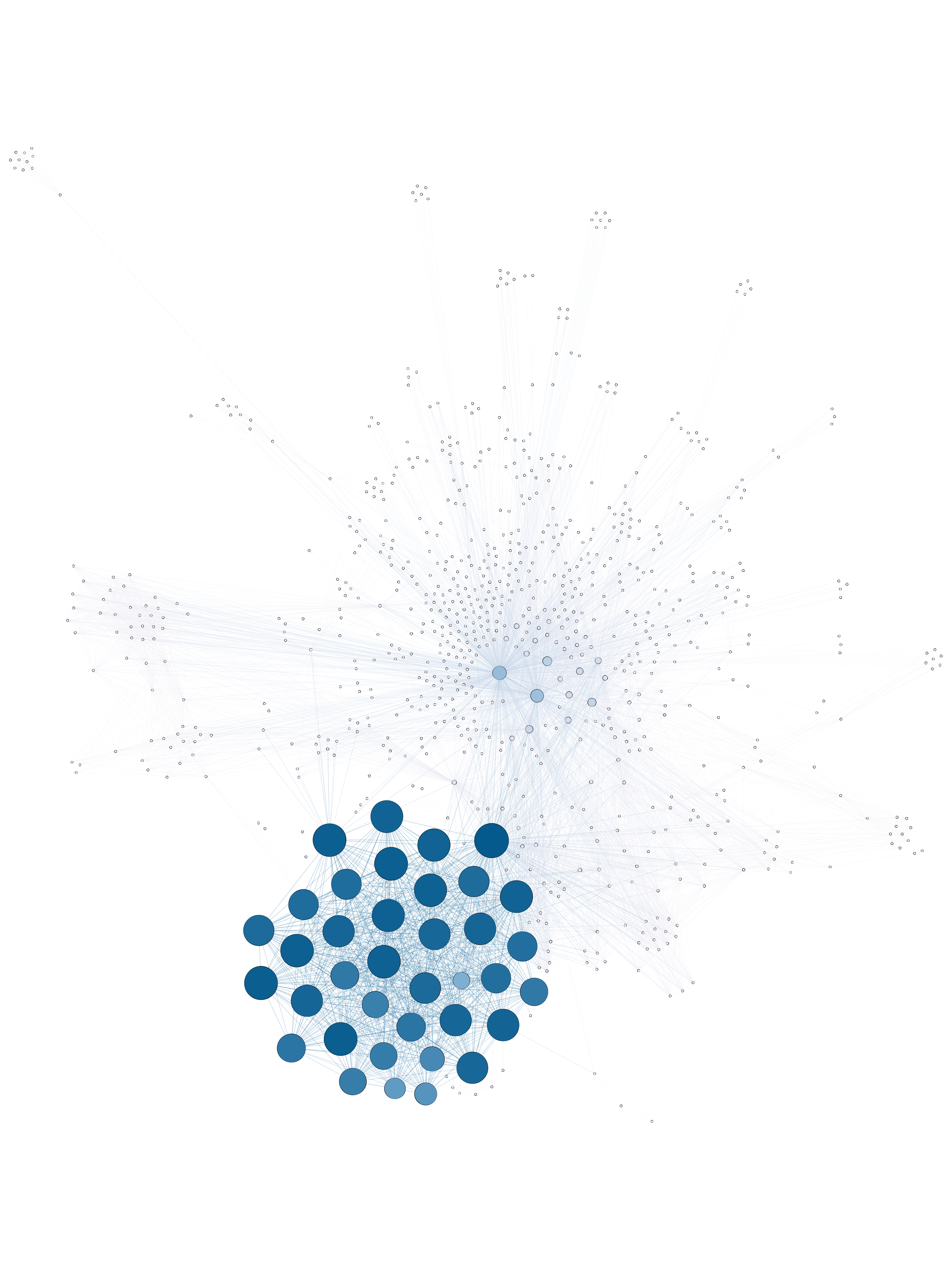}
            
            {{\small NBT Katz centrality $t = 0.95/\rho({B^{\circ 1/2}})$}}    
            \label{fig:NBT95}
        \end{subfigure}
        \caption[]
        {\small Visualizations of classical \newnew{(top/red)} and NBT Katz \newnew{(bottom/blue)} across the static email network with large node size and dark colour indicating large centrality values; darker edges indicate a larger weight.} 
        \label{fig:NetVis}
    \end{figure*}

\begin{figure}[t]
    \centering
    \begin{subfigure}{\textwidth}
        \includegraphics[width=\textwidth]{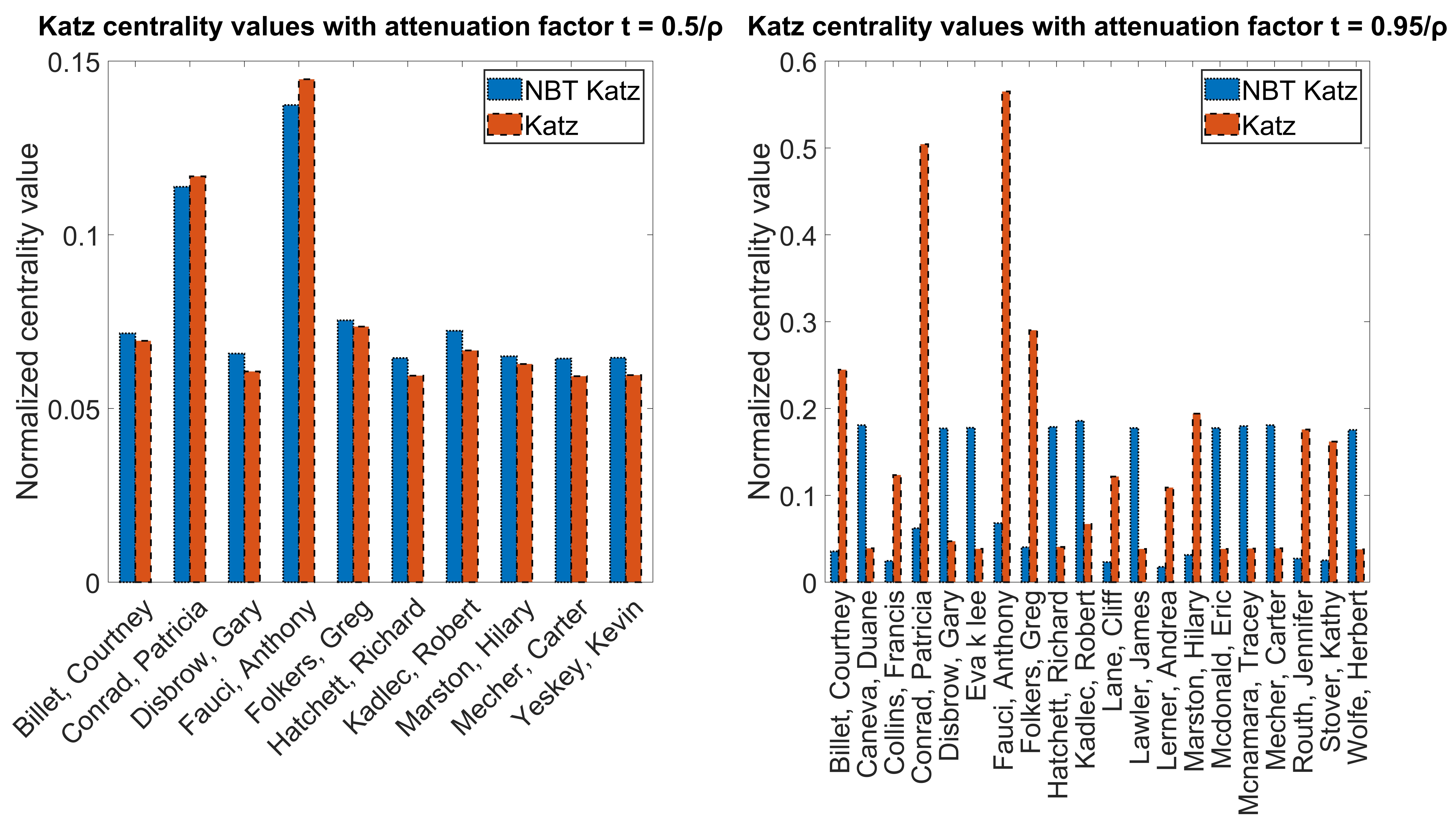}
    \end{subfigure}
    \caption{Classical and nonbacktracking Katz centrality vector values for backtracking fully forbidden with attenuation factor $t =  0.5/{\rho}$ and $t= 0.95/{\rho}$, respectively. In each plot we display the union of the 10 most central nodes according to each centrality measure.}
    \label{fig:Barstatic}
\end{figure}

\begin{figure}[ht]
    \centering
    \begin{subfigure}{\textwidth}
    \includegraphics[width=\textwidth]{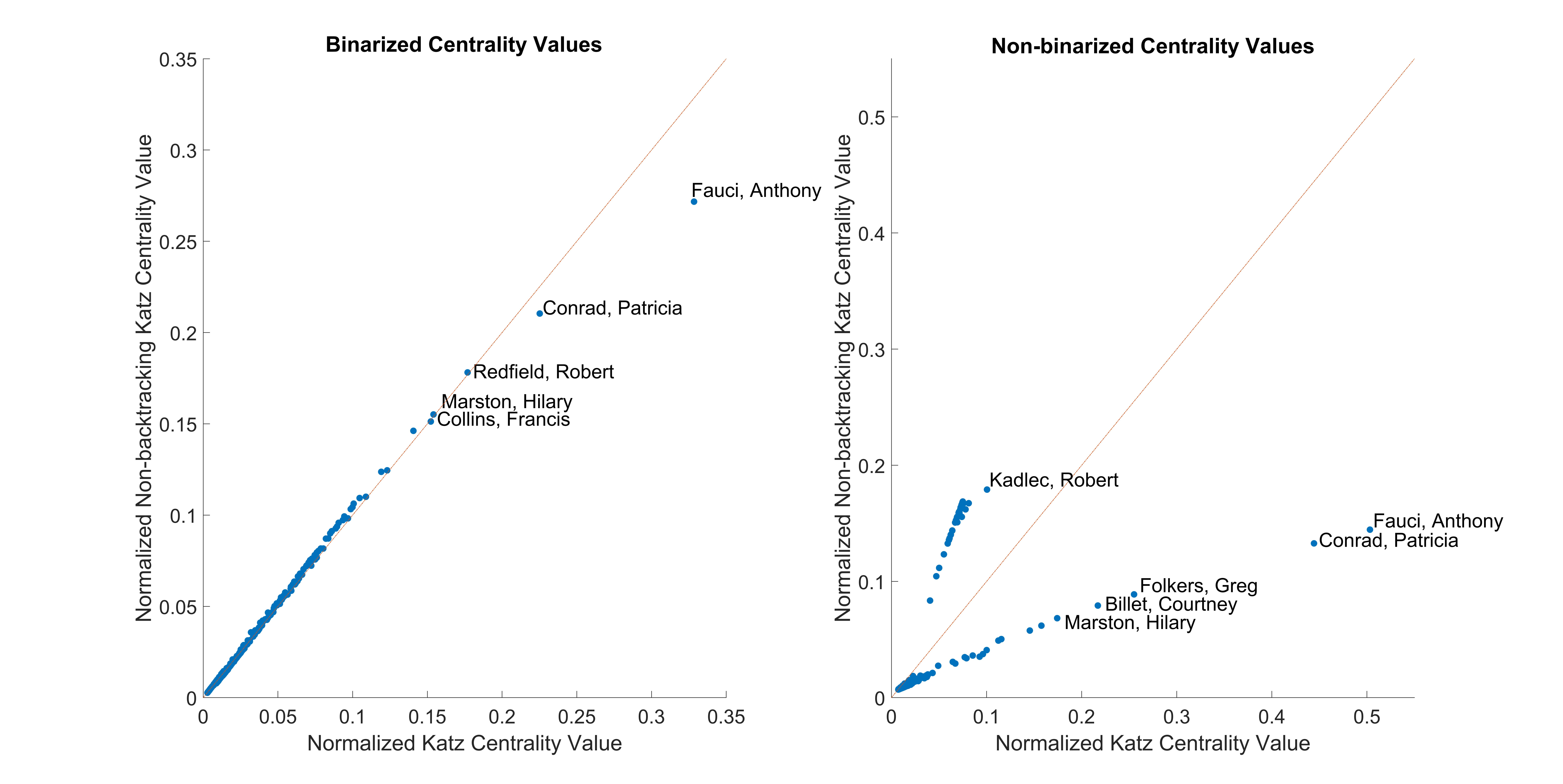}
    \end{subfigure}
    \caption{Scatter plots of normalized NBT Katz centrality against normalized classical Katz centrality corresponding to 
    binarized and non-binarized static networks with attenuation factor $t=0.95/\rho$.}
    \label{fig:binscatter}
\end{figure}

The results are visualised in 
Figures~\ref{fig:NetVis}, \ref{fig:Barstatic} and \ref{fig:binscatter}. Figure~\ref{fig:NetVis} shows that NBT Katz centrality 
emphasizes a clique not containing the node corresponding to Antony Fauci,
and that for large values of the attenuation factor this clique begins to dominate the ranking to such an extent that the node corresponding to Anthony Fauci, which occupies the central position in the network visualisation, is no longer counted among the 10 most central nodes. This can be seen in Figure~\ref{fig:Barstatic} which depicts the nonbacktracking and classical Katz normalized centrality values for the union of the 10 most central nodes in the static network. The left bar chart in Figure~\ref{fig:Barstatic} indicates that both classical and nonbacktracking Katz agree on the 10 most central nodes of which `Anthony Fauci' is most central when $t = 0.5/\rho$. However the rightmost figure depicts a complete divergence in the ten most highly-ranked nodes produced by classic and nonbacktracking Katz centralities respectively. In particular we see that while the `Anthony Fauci' node remains fairly central according to both measures, nodes belonging to the clique shown in Figure~\ref{fig:NetVis} have overtaken it in the ranking induced by nonbacktracking Katz centrality. The clique identified in this case consists exclusively of participants (i.e., either directly sent or received an email within the thread, or were CC'd in an email within the thread) in the so-called `Red Dawn' email thread that was used throughout the pandemic ``to provide thoughts, concerns, raise issues, share information across various colleagues responding to Covid-19'' \cite{Leopold-2021-fauci-emails}.

The effect of weighted edges on the rankings produced by nonbacktracking and classical Katz centralities for the static network is demonstrated in Figure~\ref{fig:binscatter}. The figure contains two scatter graphs of the normalized nonbacktracking Katz centrality vector ($t = 0.95/\rho(B^{\circ 1/2})$) plotted against the Katz centrality vector ($t = 0.95/\rho(A)$) for both the original network (right) and a binarized modified network (left), which is formed from the original network by setting all edge weights to $1$.

In particular we see that the presence of non-uniformly weighted edges in the network produces greater variation in the nonbacktracking and classical Katz centrality vectors.


\subsection{Analysis on Temporal Networks}
We now move on to the case of a time-dependent network, and we note that the special case of an unweighted network with backtracking permitted corresponds to the work in \cite{grindrod2011communicability} wherein the {\it dynamic communicability matrix} $\mathcal{Q}(t)$ associated to such a network is defined as the product of the successive resolvents 
\begin{equation}\label{def:Q}
    \mathcal{Q}(t) = (I-tA^{[1]})^{-1} (I-tA^{[2]})^{-1} \cdots (I-tA^{[N]})^{-1}.
\end{equation}
Here $A^{[i]}$ is the adjacency matrix associated to the $i$-th time-stamp of the temporal network $\mathcal{G}$. Katz centrality can then be computed via the formula
\begin{equation}\label{def:BTKatz}
   \boldsymbol{x}(t) = \mathcal{Q}(t)\boldsymbol{1}.
\end{equation}
This formula accounts for all walks across the temporal network $\mathcal{G}$ including those that backtrack in space and between time-stamps.


\newnew{The temporal network $\mathcal{G}$ analyzed in this section is the largest temporal strong component \cite{bhadra2003complexity} of the provided email data, i.e., the largest component that is connected in the sense that there exists a time-respecting path between any two nodes contained within.
This network consists of a collection of 100 directed networks associated with the date 2018-09-04 and the 99 consecutive days between 2020-01-26 and 2020-05-05. In this network we have a directed weighted edge $(i,j) \in E(G^{[\tau_t]})$, if node $j$ is a recipient of, or is CC'd in, an email sent by node $i$. The weight of such an edge is equal to the number of such emails sent during the $t$-th timestamp.} 

We reiterate here that when treating temporal networks there is a range of possible nonbacktracking regimes, as outlined in Definition \ref{def:M}. The choice of appropriate backtracking regime is highly context-dependent. 
For the data set analyzed here, it is reasonable to forbid backtracking entirely, since the time-stamps associated with the temporal network have an almost uniform spacing of one day, and the time taken to reply to an email is on a similar scale to the spacing between time-stamps. It is worth mentioning that this choice to fully forbid backtracking is subjective and other regimes may also be reasonable.


Our analysis of the spectrum of the global temporal transition matrix $M$ associated to the graph $\mathcal{G}$ with backtracking fully-forbidden yields the permitted ranges of attenuation factor $t$ shown in Table~\ref{tab:NBTconv}. We contrast this with the permitted range of $t$ in the case of classical Katz centrality via the dynamic communicability matrix $\mathcal{Q}$ as defined in \eqref{def:Q}.

\begin{table}[t]
    \centering
        \caption{Temporal network convergence information.}
    \begin{tabular}{c|c}
    \hline 
        $\rho(M)$ &    $5.025$  \\ 
         \text{nonbacktracking permitted range of $t$} & $t \in [0,0.1990)$   \\ 
         \hline 
         \hline
         $\max_i(\rho(A^{[i]}))$ &  $8.832$\\ 
         \text{Backtracking permitted range of $t$} &  $t \in [0,0.1132)$ \\
         \hline
    \end{tabular}
    \label{tab:NBTconv}
\end{table}

\begin{figure}[t]
    \centering
    \begin{subfigure}[b]{0.495\textwidth}
        \includegraphics[width=\textwidth]{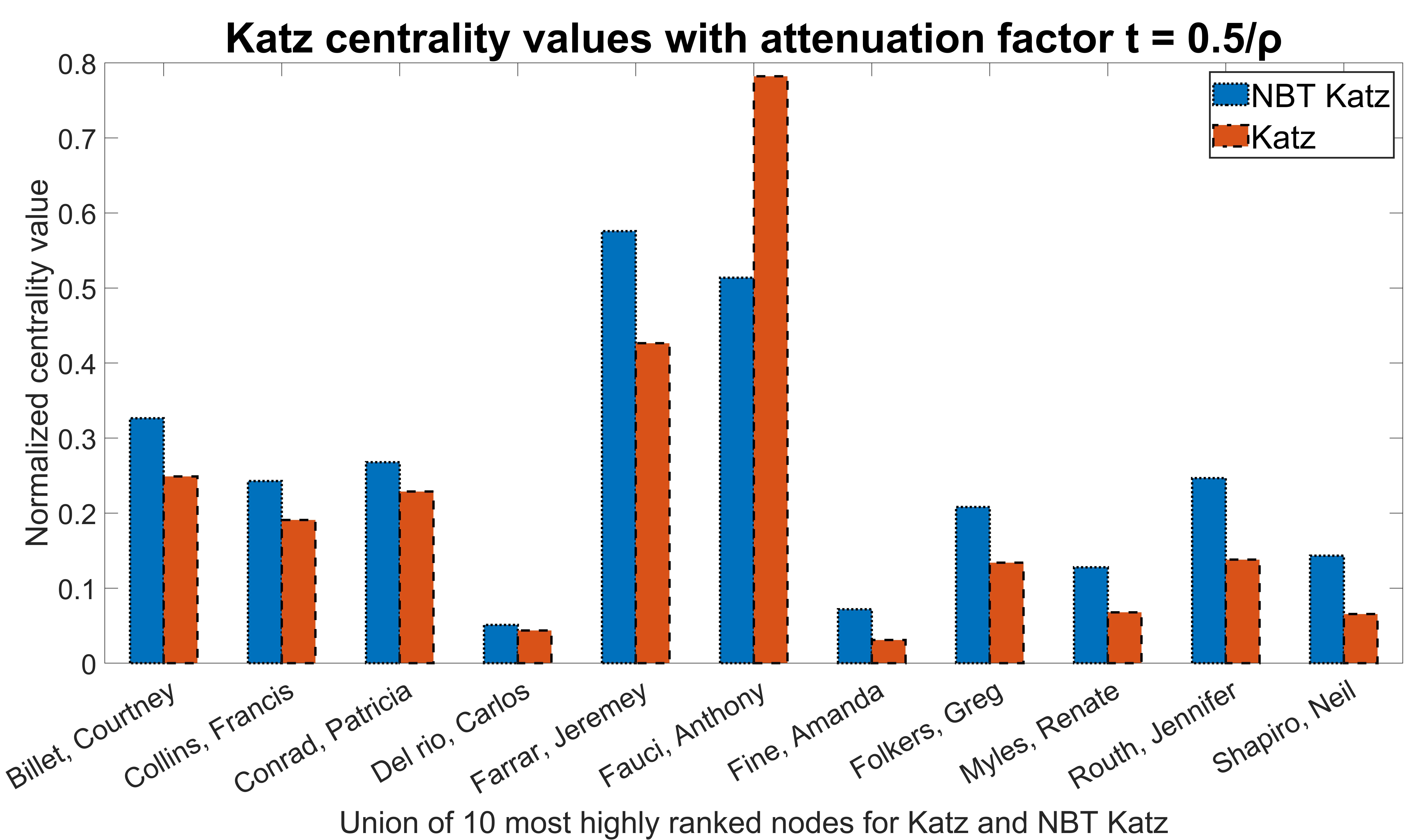}
    \end{subfigure}
        \begin{subfigure}[b]{0.495\textwidth}
        \includegraphics[width=\textwidth]{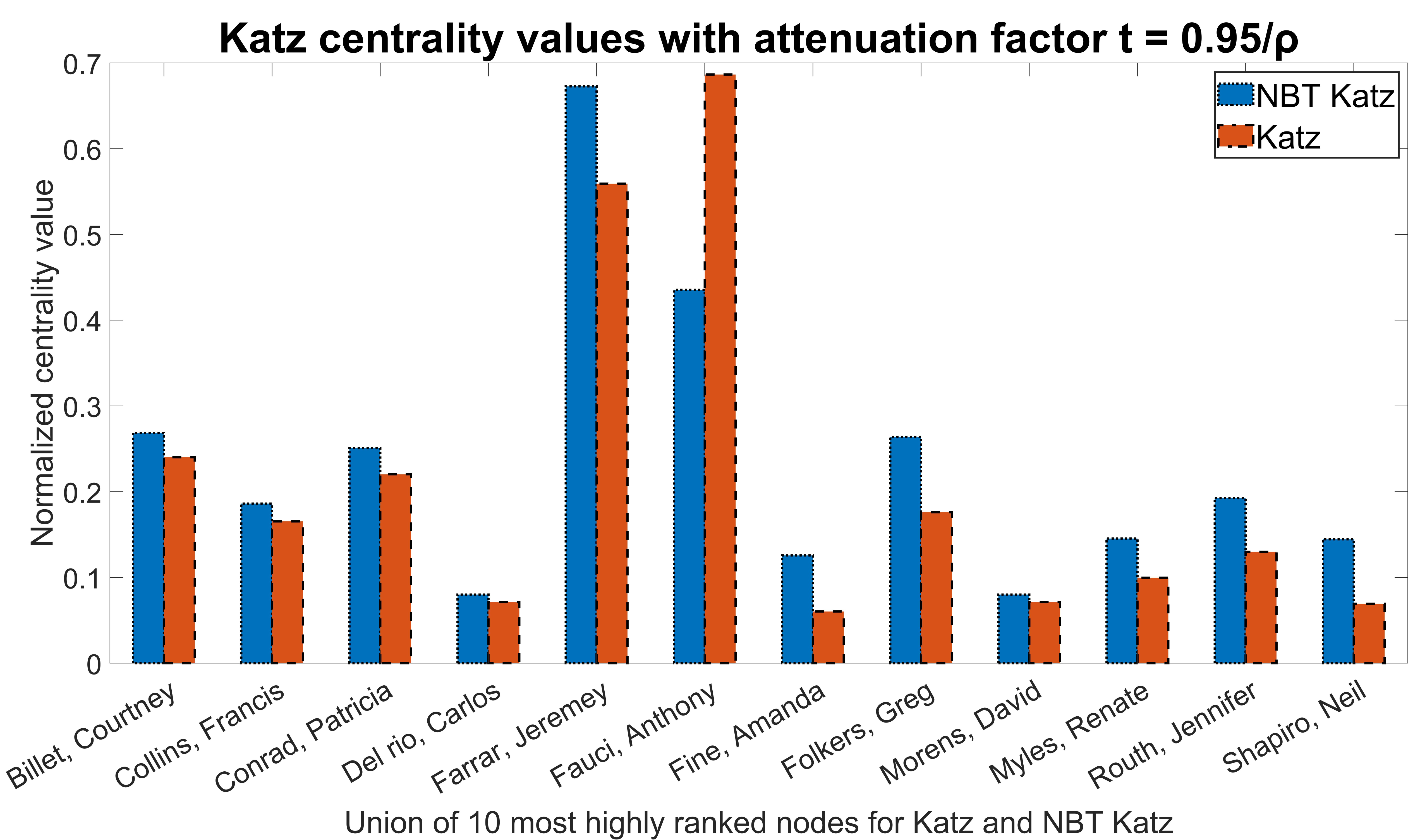}
    \end{subfigure}
        \caption{The time-evolving network centrality vector values for both NBT and classical Katz with attenuation factor $t =  0.5/\rho$ and $t= 0.95/{\rho}$ respectively. In each plot we display the union of the ten most central nodes according to each centrality measure.}
        \label{fig:NBTtemporalbar}
\end{figure}

Figure~\ref{fig:NBTtemporalbar} depicts two bar charts which display the normalized centrality values for both classical and nonbacktracking Katz centralities for $t = 0.5/\rho$ and $t = 0.95/\rho$ respectively, where $1/\rho$ is the upper-limit of the respective regime as given in Table~\ref{tab:NBTconv}. In particular $1/\rho$ is equal to $1/\rho(M)$ \newnew{(see the proof of \cite[Theorem 5.2]{manmevanni})} in the case of nonbacktracking Katz centrality, where $M$ is the matrix described in Definition~\ref{def:M} {\it (iv)} that is, the form of $M$ in which all forms of backtracking are forbidden. In the case of classical Katz centrality $1/\rho$ is given by $1/\max_i(\rho(A^{[i]})$, the reciprocal of the largest principal eigenvalue of the adjacency matrices.  
In Figure~\ref{fig:NBTtemporalbar} we report results for 12 nodes, which are selected by taking the union of the 10 most highly ranked nodes for classical Katz and the 10 most highly ranked nodes for NBT Katz, when $t = 0.95/\rho$. 

\begin{figure}[hbtp]
        \centering
    \begin{subfigure}[b]{\textwidth}
\includegraphics[width=\textwidth]{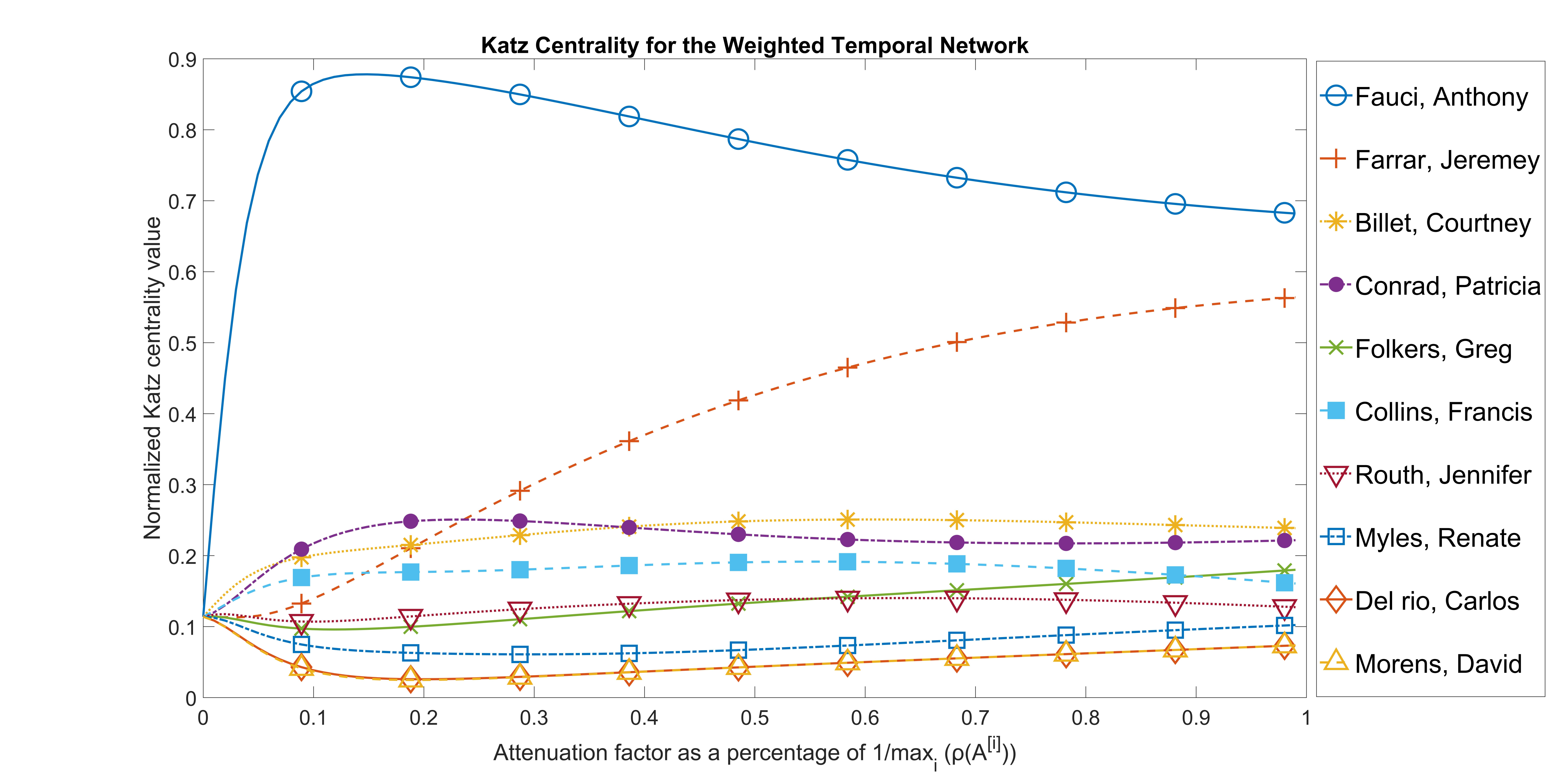}
    \end{subfigure}
    \hfill
        \begin{subfigure}[b]{\textwidth}
        \includegraphics[width=\textwidth]{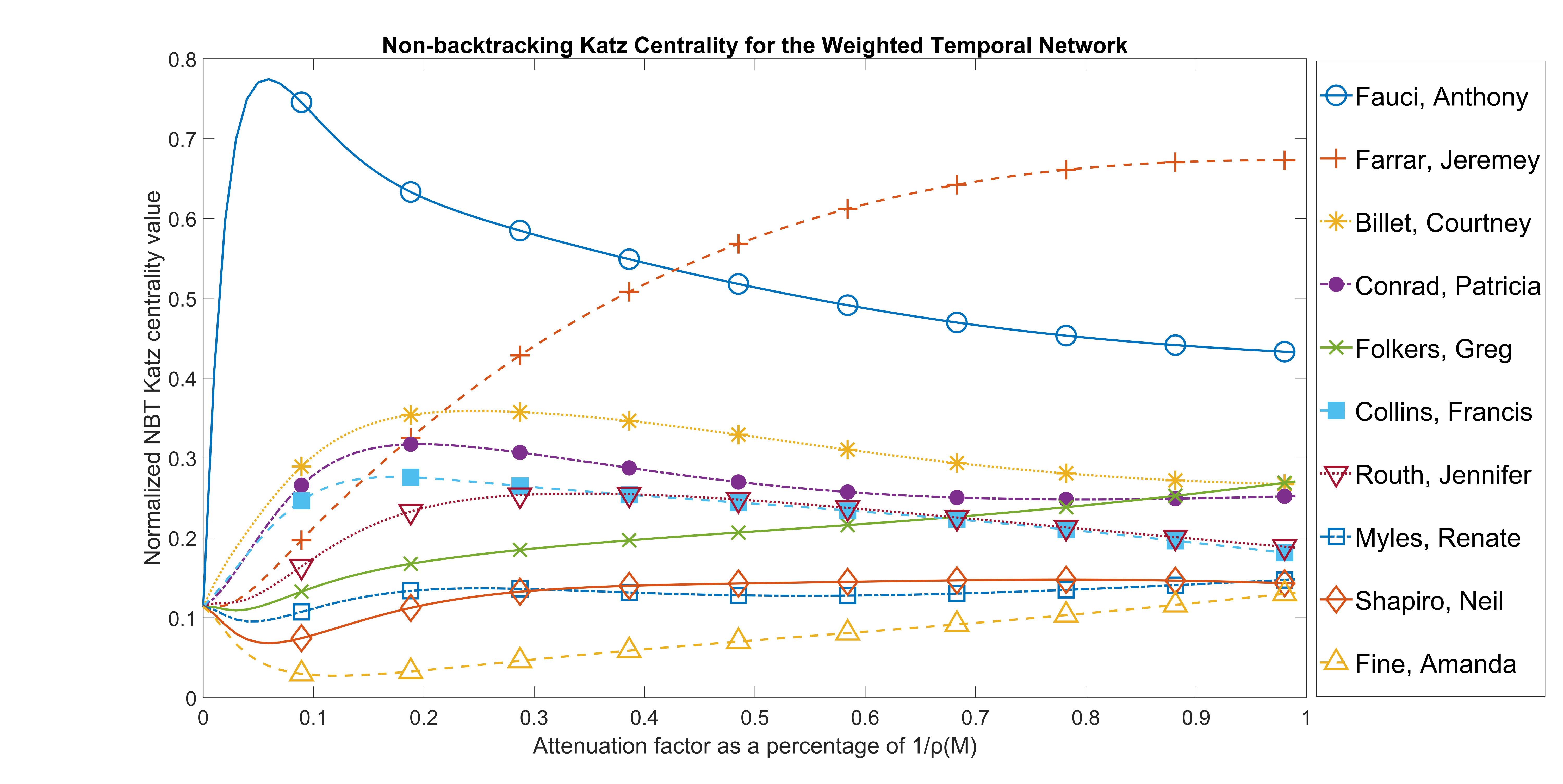}

    \end{subfigure}
        \caption{Plots of the normalized Katz (upper) and nonbacktracking Katz (lower) centralities vector values for 10 most prominent nodes  (i.e., those with the largest centrality value as of the upper limit of the attenuation factor $t$) within the weighted temporal network}
    \label{fig:contrastKweighted}
\end{figure}

In Figure~\ref{fig:contrastKweighted} we plot for the weighted temporal network both the classical and nonbacktracking Katz centrality values of 10 selected nodes against the attenuation factor $t$ which ranges from $0\%$ to $99\%$ of its permitted range (as given in Table~\ref{tab:NBTconv}). The 10 nodes were selected such that they are the most central for large values of $t$.

\begin{figure}[htbp]
    \centering
            \begin{subfigure}[b]{\textwidth}
        \includegraphics[width=\textwidth]{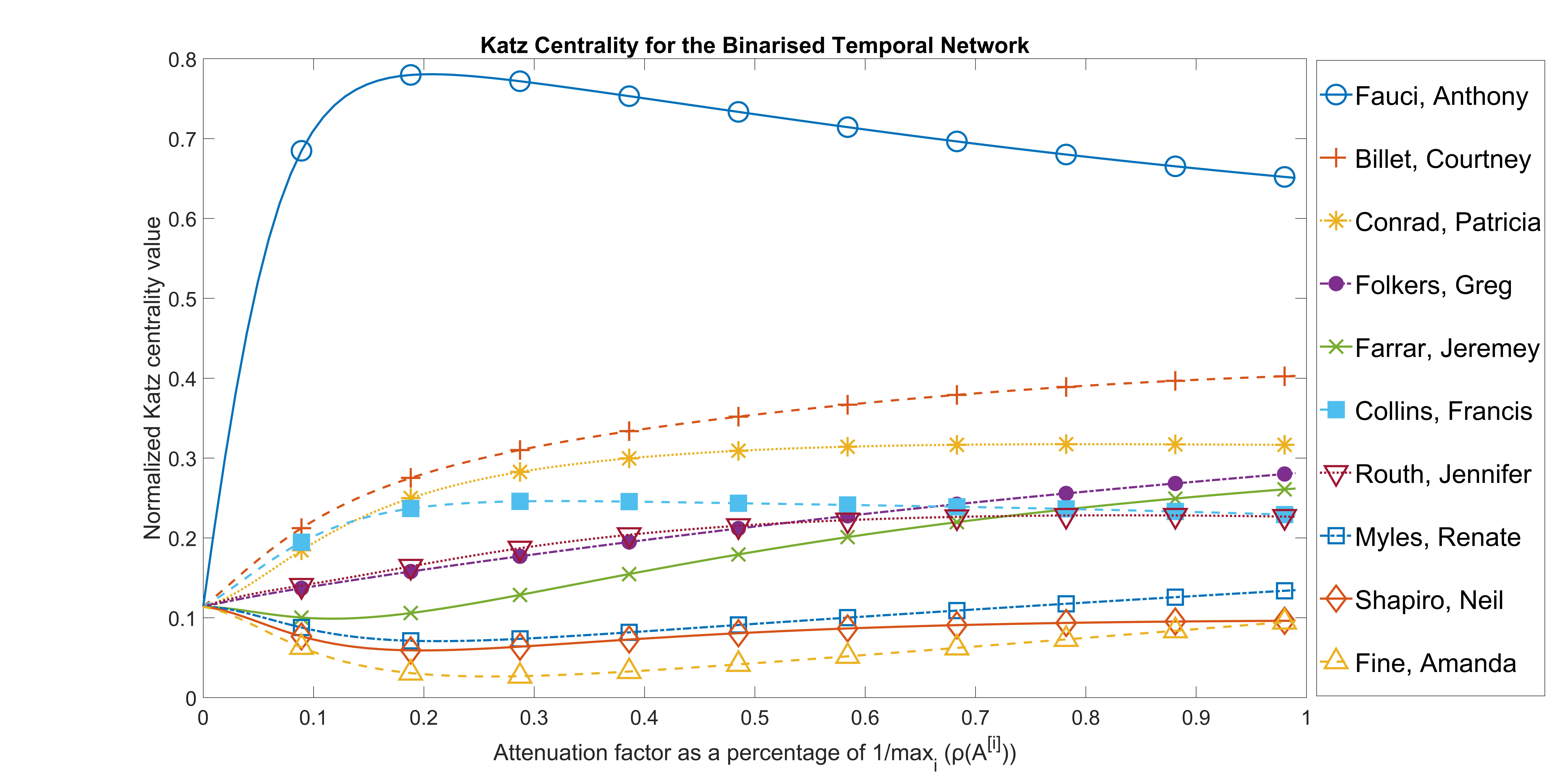}
        
        \end{subfigure}
        \hfill
        \begin{subfigure}[b]{\textwidth}
        \includegraphics[width=1\textwidth]{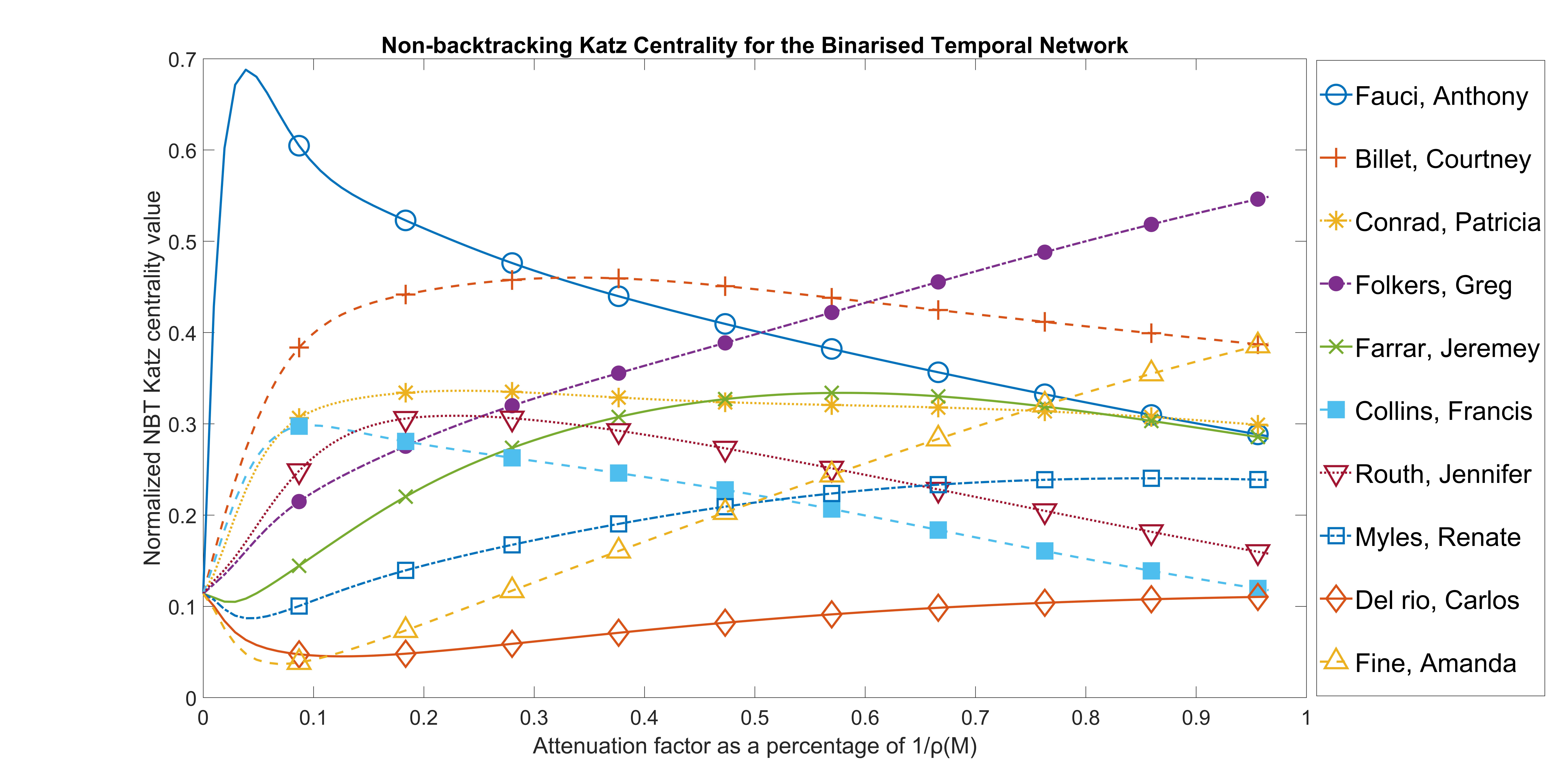}
                
        \end{subfigure}

    \caption{Plots of the normalized Katz (upper) and nonbacktracking Katz (lower) centralities vector values for 10 most prominent nodes  (i.e., those with the largest centrality value as of the upper limit of the attenuation factor $t$) within the binarized temporal network.}
    \label{fig:contrastKbinary}
\end{figure}

Figure~\ref{fig:contrastKbinary} presents results for the same experiment, this time carried out with the binarized version of the temporal network, i.e., the temporal network with all non-zero weights set to $1$. 

It is interesting to note that nonbacktracking Katz identifies a node distinct from ``Anthony Fauci'' as the most central node for large values of $t$, favouring instead the node ``Jeremy Farrar'' which is considerably lower ranked in the static networks produced from the same data set. Furthermore by comparing Figures~\ref{fig:contrastKweighted} and \ref{fig:contrastKbinary}, we observe the large effect that weighting has on the two centrality measures.

\section{Discussion}
\label{sec:disc}
Our aim in this work was to develop a useful theory 
for the enumeration of nonbacktracking walks \newnew{as well as for associated centrality measures,} in the case of \newnew{edge weights that are 
combined 
multiplicatively}.
We showed in Theorem~\ref{thm:recurrence} that in contrast to 
the unweighted case where a four-term recurrence is sufficient 
to count nonbacktracking walks of different lengths, the 
weighted case gives rise to a recurrence where the walk count at length $k$ depends on walk counts for all shorter lengths.
Despite this added complexity, the 
resulting formulas for the 
standard generating function
in 
Theorem~\ref{thm:genkatz}
and corresponding node centrality measure 
in
Corollary~\ref{cor:cent} 
are straightforward to evaluate.

We also showed in 
Theorem~\ref{thm2} 
that when working at the line graph level, the introduction of 
appropriate componentwise square roots allows us to 
develop a theory that extends to the unweighted case, with Theorem~\ref{thm:line}
summarizing the results, and 
Theorem~\ref{thm:fcenttime} 
dealing with more general 
time-evolving graph sequences.

 A practical take-home message is that
 a theory of
 nonbacktracking walk counts 
 for static or dynamic weighted graphs is available, 
 with 
 corresponding 
 computational algorithms that have the same complexity 
 as in the unweighted case.

\section*{Acknowledgements}
We thank the Editor and the Referees for their useful comments.

\printbibliography

@article{JK,































        AUTHOR="Jon Kleinberg",































        TITLE="Navigation in a small world",































        JOURNAL = "Nature",































        VOLUME = "406",































        YEAR="2000",































        PAGES="845"































}

@article{Katz53,















 author = "L. Katz",















 year = "1953",















 title = "A new index derived from sociometric data analysis",















 journal = "Psychometrika", 















 volume = "18",















 pages = "39--43"















}

@article{EV05,
author = "E. Estrada and J. A. Rodr{\'i}guez-Vel{\'a}zquez",
title = "Subgraph centrality in complex networks",
journal = "Physical Review E",
year = "2005",
volume = "71",
pages = "056103"
}

@book{Newmanbook,
author = "M. E. J. Newman",
title = "Networks: an Introduction",
publisher = "Oxford Univerity Press",
address = "Oxford",
year = "2010"
}

@ARTICLE{holme11,



   author = {{Holme}, P. and {Saram{\"a}ki}, J.},



   title = "Temporal Networks",



   journal = "Physics Reports",



   year = "2012",
  
 volume = "519",

  pages = "97--125"


}

@article{MS01,
author={H. Mizuno and I. Sato},
title={Zeta functions of digraphs},
journal={Linear Algebra and its Applications},
year={2001},
volume={336},
number={1--3},
pages={181--190}
}

@article{KKS21,
author={T. Komatsu and N. Konno and I. Sato},
title={A zeta function with respect to non-backtracking alternating walks for a digraph},
journal={Linear Algebra and its Applications},
year={2021},
volume={620},
pages={344--367}
}

@book{Estradabook,

author = "E. Estrada",

title = "The Structure of Complex Networks",

publisher = "Oxford University Press",

address = "Oxford",

year = "2011"

}

@article{BK13,
  title={Total communicability as a centrality measure},
  author={Benzi, Michele and Klymko, Christine},
  journal={Journal of Complex Networks},
  volume={1},
  number={2},
  pages={124--149},
  year={2013},
  publisher={Oxford University Press}
}

@misc{manmevanni,
      title={Generating functions of non-backtracking walks on weighted digraphs: radius of convergence and Ihara's theorem}, 
      author={Vanni Noferini and Mar\'{i}a C. Quintana},
      year={2023},
      eprint={https://arxiv.org/pdf/2307.14200.pdf},
      archivePrefix={arXiv},
      primaryClass={math.CO}
}

@article{Kempton16,
author = "Kempton, Mark",
title = "Non-Backtracking Random Walks and a Weighted Ihara’s Theorem,",
journal = "Open J. Discrete Math.",
volume = "6",
year = "2016", 
pages = "207--226"
}

@article{ST96,
title = "Zeta functions of finite graphs and coverings",
journal = "Advances in Mathematics",
volume = "121",
number = "1",
pages = "124--165",
year = "1996",
author = "H.M. Stark and A.A. Terras"
}

@article{TP09,
author = "Tarfulea, Andrei and Perlis, Robert",
title = "An {I}hara formula for partially directed graphs",
journal = "Linear Algebra and its Applications",
volume = "431",
year = "2009",
pages = "73--85"
}

@incollection{HST,
   author = {Horton, Matthew D. and Stark, H. M. and Terras, Audrey A.},
   title = {What are zeta functions of graphs and what are they good for?},
   booktitle = {Quantum graphs and their applications},
   series = {Contemp. Math.},
   editor = {Berkolaiko, G. and Carlson, R. and Fulling, S. A. and Kuchment, P.},
   pages = {173-190},
   volume = {415},
   year = {2006}
}

@inproceedings{BL70,
        author = "Bowen, R. and Lanford, O. E.",
        title = "Zeta functions of restrictions of the shift transformation",
 booktitle = "Global Analysis: Proceedings of the Symposium in Pure Mathematics of the Americal Mathematical Society, University of California, Berkely, 1968",
        year = "1970",
 publisher = "American Mathematical Society",
editor = "Shiing-Shen Chern and Stephen Smale",
pages = "43--49"
}

@article{AHN19b,
	title = "Non-backtracking {P}age{R}ank",
	author = "Arrigo, Francesca and Higham, Desmond J. and Noferini, Vanni",
	journal = "Journal of Scientific Computing",
	year = "2019",
	doi="https://doi.org/10.1007/s10915-019-00981-8"
}

@article{AGHN17a,
	title = "Non-backtracking walk centrality for directed networks",
	author = "Arrigo, Francesca and Grindrod, Peter and Higham, Desmond J. and Noferini, Vanni",
	journal = "Journal of Complex Networks",
	year = "2018",
	volume = "6",
	number = "1",
	pages = {54--78},
	publisher = "Oxford University Press"
}

@article{AGHN17b,
	title = "On the exponential generating function for non-backtracking walks",
	author = "Arrigo, Francesca and Grindrod, Peter and Higham, Desmond J. and Noferini, Vanni",
	journal = "Linear Algebra and Its Applications ",
	year = "2018",
	volume = "79",
	number = "3",
	pages = "781--801"
}

@article{GHN18,
	title = "The deformed graph {L}aplacian and its applications to network centrality analysis",
	author = "Grindrod, Peter and Higham, Desmond J. and Noferini, Vanni",
	year = "2018",
	journal = "SIAM Journal on Matrix Analysis and Applications",
	volume = "39",
	number = "1", 
	pages = "310--341",
	publisher = "SIAM"
}

@article{Hash90,
author = "Hashimoto, K.",
year = "1990",
title = "On {Z}eta and {L}-functions of finite graphs",
journal = "Intl. J. Math.", 
volume = "1",
pages = "381--396"
}

@article{MRMPO10,
title = "Community structure in time-dependent, multiscale, and multiplex networks",
author = "P. Mucha and T. Richardson and K. Macon and M. Porter and J. Onnela",
journal = "Science",
volume = "328",
pages = "876--878",
year = "2010"
}

@article{fenu2017block,
  title={Block matrix formulations for evolving networks},
  author={Fenu, Caterina and Higham, Desmond J},
  journal={SIAM Journal on Matrix Analysis and Applications},
  volume={38},
  number={2},
  pages={343--360},
  year={2017},
  publisher={SIAM}
}

@article{TSE19a,
title = "Non-backtracking cycles: {L}ength spectrum theory and graph mining applications",
author = "Leo Torres and Pablo Suarez-Serrato and Tina Eliassi-Rad",
journal = "Journal of Applied Network Science",
volume = "41", 
year = "2019"
}

@article{TCTE21,
author = "L. Torres and K. S. Chan and H. Tong and T. Eliassi-Rad",
title = "Nonbacktracking eigenvalues under node removal: {X}-centrality and targeted immunization",
journal = "SIAM Journal on Mathematics of Data Science",
volume = "3",
year = "2021"
}

@article{TCDM21,
  title = {Approximating nonbacktracking centrality and localization phenomena in large networks},
  author = {Tim\'ar, G. and da Costa, R. A. and Dorogovtsev, S. N. and Mendes, J. F. F.},
  journal = {Phys. Rev. E},
  volume = {104},
  pages = {054306},
  year = {2021}
}

@article{VDF09,
title = "Similarity matrices for colored graphs",
author = "Paul {Van Doreen} and Catherine Fraikin",
journal = "Bull. Belg. Math. Soc. Simon Stevin",
volume = "16",
year = "2009",
pages = "705--722"
}

@article{AHNW22,
  title={Dynamic Katz and related network measures},
  author={Arrigo, Francesca and Higham, Desmond J and Noferini, Vanni and Wood, Ryan},
  journal={Linear Algebra and its Applications},
  year={2022},
  publisher={Elsevier}
}

@article{Beyond,
title = "Beyond non-backtracking: non-cycling network centrality measures",
author = "Arrigo, Francesca and Higham, Desmond J. and Noferini, Vanni",
journal = "Proc. R. Soc. A",
volume= "476",
number = "2235",
pages = "20190653",
year = "2020"
}

@misc{benson2021fauciemail,
      title={fauci-email: a json digest of Anthony Fauci's released emails}, 
      author={Austin R. Benson and Nate Veldt and David F. Gleich},
      year={2021},
      eprint={2108.01239},
      archivePrefix={arXiv},
      primaryClass={cs.SI}
}

@article{grindrod2011communicability,
  title={Communicability across evolving networks},
  author={Grindrod, Peter and Parsons, Mark C and Higham, Desmond J and Estrada, Ernesto},
  journal={Physical Review E},
  volume={83},
  number={4},
  pages={046120},
  year={2011},
  publisher={APS}
}

@inproceedings{bhadra2003complexity,
  title={Complexity of connected components in evolving graphs and the computation of multicast trees in dynamic networks},
  author={Bhadra, Sandeep and Ferreira, Afonso},
  booktitle={International Conference on Ad-Hoc Networks and Wireless},
  pages={259--270},
  year={2003},
  organization={Springer}
}

@misc{Leopold-2021-fauci-emails,
  title = {Anthony Fauci’s Emails Reveal The Pressure That Fell On One Man},
  author = {Natalie Bettendorf and Jason Leopold},
  % howpublished = {BuzzFeed News, \url{https://www.buzzfeednews.com/article/nataliebettendorf/fauci-emails-covid-response}},
  month = {June},
  year = {2021},
  url = {https://s3.documentcloud.org/documents/20793561/leopold-nih-foia-anthony-fauci-emails.pdf}
}

\end{document}